\documentclass[10pt,reqno]{amsart}
\usepackage{amsthm, amscd, amssymb, amsfonts, amsmath, dsfont, mathrsfs}
\usepackage[english]{babel}
\usepackage[latin1]{inputenc}
\usepackage{graphicx}
\usepackage{enumerate}
\usepackage{color}
\usepackage[usenames]{xcolor}
\usepackage{hyperref}
\usepackage{multicol}

\newtheorem{thm}{Theorem}
\newtheorem{cor}{Corollary}
\newtheorem{lem}{Lemma}
\newtheorem{prop}{Proposition}
\newtheorem{rem}{Remark}

\newcommand{\R}{\mathds{R}}
\newcommand{\p}{p(\cdot)}
\newcommand{\q}{q(\cdot)}
\newcommand{\pri}{p'(\cdot)}
\newcommand{\loc}{\text{\upshape{loc}}}
\newcommand{\vertiii}[1]{{\left\vert\kern-0.25ex\left\vert\kern-0.25ex\left\vert #1 
		\right\vert\kern-0.25ex\right\vert\kern-0.25ex\right\vert}}

\title[Boundedness of fractional operators associated ...]{Boundedness of fractional operators associated with Schr\"odinger operators on weighted variable Lebesgue spaces via extrapolation}

\author{R. Ayala and A. Cabral}



\thanks{This research is partially supported by grants from Consejo Nacional de Investigaciones Cient\'ificas y T\'ecnicas (CONICET), Universidad Nacional del Litoral (UNL) and Universidad Nacional del Nordeste (UNNE), Argentina.}

\address{Departamento de Matem\'atica, Facultad de Ingenier\'ia Qu\'imica, CONICET, UNL, 3000 Santa Fe, Argentina.}
\email{rocioayalazara@gmail.com}

\address{Instituto de Modelado e Innovaci\'on Tecnol\'ogica, CONICET-UNNE, and Departamento de Matem\'atica, Facultad de Ciencias Exactas y Naturales y Agrimensura, UNNE, 3400, Corrientes, Argentina.}
\email{enrique.cabral@comunidad.unne.edu.ar}

\date{}

\begin{document}
	
\begin{abstract}
In this work we obtain boundedness results for fractional operators associated with Schr\"odinger operators $\mathcal{L}=-\Delta+V$ on weighted variable Lebesgue spaces. These operators include fractional integrals and their respective commutators. Particularly, we obtain weighted inequalities of the type $L^{\p}$--$L^{\q}$ and estimates of the type $L^{\p}$--Lipschitz variable integral spaces. For this purpose, we developed extrapolation results that allow us to obtain boundedness results of the type described above in the variable setting  by starting from analogous inequalities in the classical context. Such extrapolation results generalize what was done in~\cite{Pola} and~\cite{BCH} for the classic case, that is, $V\equiv0$ and $\p$ constant, respectively.

\end{abstract}

\maketitle

\section{Introduction}

Let us consider the second-order Schr\"odinger differential operator in $\R^n$ with $n\ge3$, defined by
\begin{equation}
\mathcal{L}=-\Delta+V,
\end{equation}
where $V\ge0$ and belongs to a reverse-H\"older class $RH_\nu$ for some exponent $\nu\ge n/2$, i.e., there exists a constant $C$ such that
 \begin{equation*}
\left(\frac{1}{|B|}\int_B V(x)^\nu dx\right)^{1/\nu}\leq \frac{C}{|B|}\int_B V(x)\,dx,
\end{equation*}						
for every ball $B\subset\R^n$.

In recent years, a wide range of operators associated to $\mathcal{L}$ have been catching the attention of several authors (see \cite{BCH, BCH-Hardy, BHS, BHS-Conm, DZ-H2, DGMTZ, Shen, Tang}). One of those is the fractional integral, which is defined by means of negative fractional powers of $\mathcal{L}$, i.e.
\begin{equation*}
\mathcal{I}_\alpha=\mathcal{L}^{-\alpha/2},
\end{equation*}
for $0<\alpha<n$.

In~\cite{DGMTZ} the authors proved that $\mathcal{I}_\alpha$ is bounded from the Lebesgue space $L^p$ to $L^q$, provided that $1/q=1/p-n/\alpha$ and $1<p<n/\alpha$, as the classical fractional integral. The corresponding weighted inequality was demonstrated in~\cite{BHS}. Since the negative powers of $\mathcal{L}$ can be expressed in terms of the heat diffusion semigroup generated by $\mathcal{L}$, the kernel of $\mathcal{I}_\alpha$ has better behavior far away from the diagonal than the classical one. Thus, the family of weights $A_{p,q}^\rho$ (see Section 3) obtained in the boundedness $L^{p}(w^p)$--$L^{q}(w^q)$ is wider than the class of weights of Muckenhoupt $A_{p,q}$  associated to the classical operator.

For the limiting case $p=n/\alpha$,  in~\cite{BCH} they characterize the weights $w$ such that  $\mathcal{I}_\alpha$ is bounded from $L^{n/\alpha}(w^{n/\alpha})$
into $BMO_{\mathcal{L}}(w)$, an appropriate weighted version of $BMO$, the John-Nirenberg space adapted to this context.

The case of the commutator $[\mathcal{I}_\alpha,b]$ of $\mathcal{I}_\alpha$ was treated in~\cite{Tang} where it was shown that it is bounded from $L^p(w^p)$ into $L^q(w^q)$, for $1/q=1/p-n/\alpha$, $1<p<n/\alpha$ and $w\in A_{p,q}^\rho$.

In this context, maximal operators of fractional type can also be considered. Its definition is closely connected to the operator $\mathcal{L}$ through the critical radius function
\begin{equation}\label{rho Schrodinger}
\rho(x)=\sup\left\{r>0:\frac1{r^{n-2}}\int_{B(x,r)}V\leq1\right\}.
\end{equation}

In particular, the reverse--H\"older condition implies that $0<\rho(x)<\infty$ for $x\in\R^n$.  Furthermore, according to \cite[Lemma 1.4]{Shen}, if $V\in RH_{\nu}$, with $\nu>n/2$, the associated function $\rho$ verifies that there exist constants $c_\rho$, $N_\rho\geq1$ such that
\begin{equation}\label{est rho}
c_\rho^{-1}\rho(x)\bigg(1+\frac{|x-y|}{\rho(x)}\bigg)^{-N_\rho}\leq\rho(y)\leq c_\rho\,\rho(x)\bigg(1+\frac{|x-y|}{\rho(x)}\bigg)^{\frac{N_\rho}{N_\rho+1}},
\end{equation}
for every $x,y\in\R^n$.

Shen uses this function $\rho$ as an essential tool to obtain boundedness on $L^p$ of many operators associated with $\mathcal{L}$, for instance the Riesz Schr\"odinger transform $\mathcal{R}=\nabla(-\Delta+V)^{-1/2}$. The point is  that when $\mathcal{R}$ is restricted to the $\rho$-neighborhood of the diagonal $\Delta_{\rho}=\{(x,y):|x-y|<\rho(x)\}$, it behaves like the classical Riesz transform $R=\nabla(-\Delta)^{-1/2}$,  in the sense that the difference of their kernels has a locally integrable singularity at the diagonal. Usually such restriction is called the $\rho$-local part of the kernel. As for the other part, the presence of the potential $V$ guarantees a good decay scaled according to $\rho$.

This function is used to describe functional spaces and weights naturally associated with $\mathcal{L}$ that are generalizations of the case $V\equiv0$.



One of the main objectives of this work is to establish, by means of extrapolation techniques, boundedness results for fractional operators associated with Schr\"odinger operators on weighted variable Lebesgue spaces.

In reference to the theory of extrapolation, in \cite{Rubio}, Rubio de Francia  proved that the classes $A_{p}$ verify an interesting property of extrapolation. More specifically, if for some $p_0\ge1$, an operator preserves  $L^{p_{0}}(w)$ for any $w\in A_{p_{0}}$, then necessarily preserves the  $L^{p}(w)$ space for every $1<p<\infty$ and every $w\in A_{p}$. Later, in \cite{Pola} Harboure,  Mac\'ias and  Segovia, proved that the $A_{p,q}$ classes have a similar extrapolation property. Moreover, they also proved that this property is not only exclusive for the boundedness between weighted Lebesgue spaces, but also it is possible to extrapolate based on a continuity behavior of the type $L^s(w)-BMO(w)$ for some $1<s<\infty$. 

Recently in~\cite{CPR}, the authors demonstrate extrapolation results that allow us to obtain continuity properties of certain operators of the type $L^p(w^p)-L^q(w^q)$ or $L^p(w^p)-\mathscr{L}_{\tilde{\delta}}(w)$ starting with hypothesis of continuity of the type $L^s(w^s)-\mathscr{L}_{\delta}(w)$ for some related parameters. Here, $\mathscr{L}_{\delta}(w)$, $0<\delta<1$, denotes the weighted integral Lipschitz space, i.e.,  the set of locally integrable function $f$ such that
\begin{equation*}\label{a2}
\frac{\|w\chi_{B}\|_{\infty}}{|B|^{1+\frac{\delta}{n}}}\int_{B}|f(x)-f_B|\,dx\leq C
\end{equation*}
holds for every ball $B\subset \R^{n}$.

Our first objective is to extend some of the results seen in~\cite{CPR} to $\mathscr{L}_{\delta}(w)$, an appropriate weighted integral Lipschitz space adapted to the Schr\"odinger context.

We are also interested in establishing extrapolation results of the type described above in the variable exponent spaces context. In this line we developed extrapolation results starting from hypothesis which involves inequalities of the type $L^s(w^s)-\mathscr{L}_{\delta}(w)$, and obtaining weighted estimates of the type $L^{\p}-L^{\q}$ or $L^{\p}-\mathscr{L}_{\delta(\cdot)}$, where the last space is a variable version of the space $\mathscr{L}_{\delta}(w)$.

The structure of this paper is the following. In Section~\ref{Sec2} we give some preliminaries and state our main theorems. In Section~\ref{Sec3} we state and prove the auxiliary results which are important tools in order to prove the theorems stated in Section~\ref{Sec2}. Later, in Section~\ref{Sec4} we  deal with the proofs of the main results. Finally, in Section~\ref{Sec5}, we give some applications to the variable exponent spaces context by mean of extrapolation techniques.

Throughout this paper, unless otherwise indicated, we will use $C$ and $c$ to denote constants, which are not necessarily the same at each occurrence.
We will say that $A\lesssim B$ when there exists a  constant $c>0$ such that $A\le cB$ and we will write $A\simeq B$ whenever $A\lesssim B$ and $B\lesssim A$.

\section{Preliminaries and main results}\label{Sec2}
From now on, we call a critical radius function to any positive function $\rho$ that satisfies \eqref{est rho}. Clearly, if $\rho$ is such a function, so it is $\beta\rho$ for any $\beta>0$.

Let's denote by $\mathcal{B}_\rho$ the family of subcritical balls of $\R^n$, i.e., the set of balls $B(x,r)$ with $x\in\R^n$ and $r\le\rho(x)$.

Taking into account what has been done in \cite{BHS} (see also~\cite{Cabral} and~\cite{Tang}) for the case of a Schr\"odinger operator, we will consider classes of weights given in terms of a critical radius function.

By a weight $w$ we mean a locally integrable function such that $0<w(x)<\infty$  a.e. We say that the weight $w$ belongs to the $A^{\rho}_p$ class , for $1<p<\infty$, if there exist constants $\theta\ge 0$ and $C>0$ such that the inequality
 \begin{equation}\label{pesos a_p}
 \left(\frac{1}{|B|}\int_B w\,dy\right)^{1/p}\left(\frac{1}{|B|} \int_B w^{-\frac{1}{p-1}}dy\right)^{1/p'}\leq C \left(1+ \frac{r}{\rho(x)}\right)^{\theta},
 \end{equation}
holds for every ball $B=B(x,r) \subset \R^{n}$. For the case $p=1$, we will say that $w$ belongs to the class  $A^{\rho}_1$ if there are constants $\theta\ge 0$ and $C>0$ such that the inequality
\[
\frac{1}{|B|}\int_B w\,dy \leq C \Big(\inf_{x \in B} w(x)\Big) \left(1+ \frac{r}{\rho(x)}\right)^{\theta},
\]
holds for every ball $B=B(x,r) \subset \R^{n}$.

On the other hand, we will say that a weight $w$ belongs to the class $A^{\rho}_{p,q}$, for $p,q\in(1,\infty)$, if there exist constants $\theta\ge 0$ and $C>0$ such that the inequality
\begin{equation}\label{pesos a_pq}
\left(\frac{1}{|B|}\int_B w^{q}\,dy\right)^{1/q}\left(\frac{1}{|B|} \int_B w^{-p'}dy\right)^{1/p'}\leq C \left(1+ \frac{r}{\rho(x)}\right)^{\theta},
\end{equation}
holds for every ball $B=B(x,r) \subset \R^{n}$. In the limiting case $q=\infty$, we will say that $w\in A^{\rho}_{p,\infty} $, if there exist $q>0$ and $C>0$ such that the inequality 
\begin{equation*}
\|w\chi_{B}\|_\infty \left(\frac{1}{|B|} \int_{B}w^{-p'}\,dy \right)^{1/p'} 
\leq C \left( 1 + \frac{r}{\rho(x)}\right)^{\theta},
\end{equation*}
holds for every ball  $B=B(x,r) \subset \R^{n}$, where $\|w\chi_B\|_{\infty}= \sup_{x \in B}w(x) $.

We will also consider two other classes of weights, namely $A_p^{\rho,\loc}$ and $A_{p,q}^{\rho,\loc}$, which are formed by those weights $w$ that satisfy~\eqref{pesos a_p} and~\eqref{pesos a_pq}, respectively, only for the balls in $\mathcal{B}_\rho$. From their definitions it is clear that $A_p^{\rho}\subset A_p^{\rho,\loc}$ and $A_{p,q}^{\rho}\subset A_{p,q}^{\rho,\loc}$,
for any $1<p<\infty$ and $1<q\le\infty$.

Some properties of these classes of weights that can be derived from their definition and that will be useful to us are the following.
\begin{prop}\label{prop: prop_w} Let $1\le p\le q\le\infty$. The following are verified 
\begin{enumerate}
\item\label{propie1} If $w\in A^{\rho}_p$, then $w \in A^{\rho}_q$.
\item\label{propie2} If  $w_1, w_2\in A^{\rho}_1$, then $w_1w_2^{1-p}\in A^{\rho}_p$.
\item\label{propie3} If $p>1$,  $w \in A^{\rho}_{p,q}$ if and only if $w^{-p'} \in A^{\rho}_{1+\frac{p'}{q}}$.
\item\label{propie4} If $p>1$ and $q<\infty$, $w \in A^{\rho}_{p,q}$ if and only if $w^{q} \in A^{\rho}_{1+\frac{q}{p'}}$.
\item\label{propie5} If $w\in A^{\rho}_{p,q}$ and $p>1$, then $w \in A^{\rho}_{p,r}$, for $q\ge r\ge p$.
\item\label{propie6} If $w\in A^{\rho}_{p,q}$ and $p>1$, then $w \in A^{\rho}_{r,q}$, for $r\ge p$.
\end{enumerate}
\end{prop}

Recall that given a weight $w$, $L^{p}(w)$ denotes the space of functions $f$ such that
\begin{equation*}
\int_{\R^{n}} |f(x)|^{p} w(x)\,dx<\infty.
\end{equation*}

The classes $A_p^\rho$ are deeply connected to the family of maximal operators $M_\rho^\theta$, defined for $\theta>0$ by means of
\begin{equation}\label{Mtheta}
M_\rho^\theta f(x)=\sup_{r>0}\bigg(1+\frac r{\rho(x)}\bigg)^{-\theta}\frac1{|B(x,r)|}\int_{B(x,r)}|f(y)|\,dy.
\end{equation}
Indeed,  these are bounded in $L^p(w)$, for $1<p<\infty$, provided that $w\in A_p^\rho$, and of weak type $(1,1)$ with respect to $w$, for $w\in A_1^\rho$ (see~Proposition 3 in~\cite{BCH} and Proposition 4.2 in~\cite{BCH-Hardy}).

As stated above we are interested in introducing Lipschitz type integral spaces in the Schr\"odinger context and establishing certain extrapolation results involving them. 

In the following, we denote by $f_B$ the average of $f$ over the $B$ ball, i.e. $\frac{1}{|B|}\int_B f$.

Given a weight $w$, we define the space $\mathscr{L}_{\delta}(w)$, for $0\leq\delta<1$, as the set of all functions $f$ on $L^{1}_{\loc}(\R^n)$ such that
\begin{equation}\label{oscilacion}
\frac{\|w\chi_B \|_\infty}{|B|^{1+\frac{\delta}{n}}} \int_B |f(y)-f_B|\,dy \leq C,\ \ \ \  \forall B \in \mathcal{B} _\rho.
\end{equation}  
and 
\begin{equation}\label{promedio}
\frac{\|w\,\chi_{B(x,\rho(x))} \|_\infty}{|B(x,\rho(x))|^{1+\frac{\delta}{n}}} \int_{B(x,\rho(x))}|f(y)|\,dy \leq C,\ \ \ \ \forall x\in\R^n.
\end{equation}

We can consider a norm in  $\mathscr{L}_{\delta}(w)$ as the smallest constant that satisfies \eqref{oscilacion} and \eqref{promedio} simultaneously, and denote it by $\vertiii{f}_{\mathscr{L}_\delta (w)}$.

\begin{rem}
A function $f$ that satisfies the inequality~\eqref{promedio} for a ball $B$, also satisfies~\eqref{oscilacion} for the same ball. This essentially follows from the fact that $\int_B|f_B|\le\int_B|f(x)|$.
\end{rem}

\begin{rem}
When $\delta=0$, the space $\mathscr{L}_{0}(w)$ coincides with the space $BMO_\rho(w)$ originally considered in~\cite{BCH}.
\end{rem}

Although our goal is to use extrapolation to prove the boundedness of some specific operators, we will state our results in a more abstract way. Following the approach established in~\cite{CUMP1} (see also~\cite{CUMP2} and~\cite{CUFMP}) we will present our extrapolation theorems for pairs of measurable, nonnegative functions $(f,g)$ belonging to some family $\mathcal{F}$.  Henceforth, if we write
\[
\|f\|_X\le\|g\|_Y,\ \ \ \ (f,g)\in\mathcal{F},
\]
where $X$ and $Y$ are spaces of functions (i.e., weighted Lebesgue spaces, classical or variable), then we mean that this inequality is true for any pair $(f,g)\in\mathcal{F}$ such that the left--hand side of this inequality is finite.

We are now in a position to state our first extrapolation results taking as hypotheses inequalities of the type $L^s(w^s)-\mathscr{L}_{\delta}(w)$. The first one gives  $L^p(w^p)-L^q(w^q)$ type estimates and the second result gives $L^p(w^p)-\mathscr{L}_{\tilde{\delta}}(w)$ estimates with $0\le\delta\le\tilde{\delta}$.

\begin{thm}\label{teo1}
Let $0\leq\delta<1$, $1<\beta<\frac{n}{\delta}$  and $\frac{1}{s}=\frac{1}{\beta}-\frac{\delta}{n}$ be such that the following inequality,
\begin{equation*}
\vertiii{f}_{\mathscr{L}_\delta(w)}\leq C \|g\,w\|_s,\ \ \ \ (f,g)\in\mathcal{F},
\end{equation*}
is verified for all $w \in A^{\rho}_{s}$ and some constant $C>0$. Then, there exists a constant $c>0$ such that the inequality
\begin{equation}\label{key1}
\|f\,w\|_q\leq c\|g\,w\|_p,\ \ \ \ (f,g)\in\mathcal{F},
\end{equation}
is satisfied for any $p$ and $q$ such that $1<p<\beta$, $\frac{1}{\beta}=\frac{1}{p}-\frac{1}{q}$ and any weight $w \in A^{\rho}_{p,q}$.
\end{thm}

\begin{thm}\label{teo2}
Let $0\leq\delta<1$, $1<\beta<\frac{n}{\delta}$ and $\frac{1}{s}=\frac{1}{\beta}-\frac{\delta}{n}$ be such that the following inequality,
\[
\vertiii{f}_{\mathscr{L}_\delta(w)}\leq C \|g\,w\|_s,\ \ \ \ (f,g)\in\mathcal{F},
\]
is verified for all $w \in A^{\rho}_{s}$ and some constant $C>0$. Then, there exists a constant $c>0$ such that the inequality
\begin{equation}\label{key2}
\vertiii{f}_{\mathscr{L}_{\tilde{\delta}}(w)}\leq c\|g\,w\|_p,\ \ \ \ (f,g)\in\mathcal{F},
\end{equation}
is true for any $0\le\tilde{\delta}\le\delta$, $\frac{1}{p}=\frac{1}{\beta}-\frac{\tilde{\delta}}{n}$ and every weight $w \in A^{\rho}_{p,\infty}$.
\end{thm}

We are also interested in obtaining extrapolation results of the type described above in the variable Lebesgue space context. We begin with some definitions and notations related to these  spaces.

Let $\p:\R^n\rightarrow[1,\infty]$ be a measurable function. Given a measurable set $A\subset\R^n$ we define
\begin{equation*}
p^{-}(A):=\text{ess}\inf_{x\in A}p(x),\hspace{1cm} p^{+}(A):=\text{ess}\sup_{x\in A}p(x).
\end{equation*}

For simplicity we denote $p^{-}=p^{-}(\R^n)$ and $p^{+}=p^{+}(\R^n)$.

 Given $\p$, the conjugate exponent $\pri$ is defined pointwise
\begin{equation*}
\frac{1}{p(x)}+\frac{1}{p'(x)}=1,
\end{equation*}
where we let $p'(x)=\infty$ if $p(x)=1$.

By $\mathcal{P}(\R^n)$ we will designate  the collection of all measurable functions $\p:\R^n\rightarrow[1,\infty]$ and by $\mathcal{P}^*(\R^n)$ the set of $p\in\mathcal{P}(\R^n)$ such that $p^{+}<\infty$.

Given  $p\in\mathcal{P}^*(\R^n)$, we say that a measurable function $f$ belongs to $L^{\p}(\R^n)$ if for some $\lambda>0$, the modular of $f/\lambda$ associated with $\p$, that is,
\begin{equation*}
\varrho_{\p}(f/\lambda)=\int_{\R^n}\left(\frac{|f(x)|}{\lambda}\right)^{p(x)}dx,
\end{equation*}
is finite. A Luxemburg type norm can be defined in $L^{\p}(\R^n)$ by taking
\begin{equation*}
\|f\|_{L^{\p}(\R^n)}=\|f\|_{\p}=\inf\{\lambda>0:\varrho_{\p}(f/\lambda)\le1\}.
\end{equation*}

These spaces are special cases of Musieliak-Orlicz spaces (see \cite{Musielak}), and generalize the classical Lebesgue spaces. For more information see, for example~\cite{CUF, Libro:D-H, Kova-Rako}.

We will denote with $L^{\p}_{\loc}(\R^n)$ the space of functions $f$ such that $f\chi_B\in L^{\p}(\R^n)$ for every ball $B\subset\R^n$.

In the case of constant exponents, the $L^p$ norm and the modular differ only by one exponent.
In variable Lebesgue spaces their relation is less direct, as shown by the following result, whose proof can be found in~\cite[Corollary~2.23]{CUF}.
\begin{prop}\label{prop: modular}
Let $\p\in\mathcal{P}^{*}(\R^n)$. Hence
\begin{enumerate}
\item if $\varrho_{\p}(f)\le C$, then $\|f\|_{\p}\le\max\{C^{1/p^{-}}, C^{1/p^{+}}\}$;
\item if $\|f\|_{\p}\le C$, then $\varrho_{\p}(f)\le\max\{C^{p^{+}}, C^{p^{-}}\}$.
\end{enumerate}
\end{prop}
\
In the classical $L^p$ spaces, $1<p<\infty$, the norm can be characterized using the identity
\begin{equation*}
\|f\|_p=\sup\int_{\R^n}f(x)g(x)\,dx,
\end{equation*}
where the supremum is considered over all functions g such that $g\in L^{p'}$ and $\|g\|_{p'}\le1$.
Analogously, we have the following result for variable Lebesgue spaces.
\begin{prop}(\cite[Theorem 2.34]{CUF})\label{prop:dualidad}
Let $\p\in\mathcal{P}(\R^n)$, $f$  a measurable function and
\begin{equation*}
|||f|||_{\p}=\sup\left\{\int_{\R^n}f(x)g(x)\,dx:\ \|g\|_{\pri}\le1\right\}.
\end{equation*}
Then,
\begin{equation*}
c\,|||f|||_{\p}\le\|f\|_{\p}\le C|||f|||_{\p},
\end{equation*}
where the constants $c$ and $C$ depend only on $\p$.
\end{prop}
On the other hand, analogously to the previous case, H\"older's inequality is also valid for variable exponents but with a constant on the right-hand side of it.
\begin{prop}\label{prop:Holder}(\cite[Lema~3.2.20]{Libro:D-H})
Given $s(\cdot), \p, \q \in\mathcal{P}(\R^n)$, be such that $1/s(\cdot)=1/\p+1/\q$. Then, for $f\in L^{\p}(\R^n)$ and $g\in L^{\q}(\R^n)$
\begin{equation*}
\|fg\|_{s(\cdot)}\le 2\|f\|_{\p}\|g\|_{\q}.
\end{equation*}
Moreover, if $s(\cdot)\equiv1$, the inequality above gives
\begin{equation*}
\int_{\R^n}|f(x)g(x)|\,dx\le 2\|f\|_{\p}\|g\|_{\pri}.
\end{equation*}
\end{prop}
\
Another elementary but useful property of the classical Lebesgue norm is that it is homogeneous in the exponent, more precisely $\|f^s\|_p=\|f\|^s_{sp}$ for $1<s<\infty$ and non-negative $f$.
This property also extends to variable Lebesgue spaces as follows.
\begin{prop}(\cite[Proposition~2.18]{CUF})\label{prop:dilatacion}
Let $\p\in\mathcal{P}^{*}(\R^n)$, so for all $s$, $1/p^{-}\le s<\infty$,
\begin{equation*}
\||f|^s\|_{\p}=\|f\|^s_{s\p}.
\end{equation*}
\end{prop}

The following conditions on the exponent arise in connection with the boundedness of the Hardy--Littlewood maximal operator $M$ in $L^{\p}(\R^n)$ (see, for example, \cite{Diening}, \cite{CUF} or~\cite{Libro:D-H}). We will say that $p$ is log--H\"older continuous, and we will write $p\in\mathcal{P}^{\log}(\R^n)$, if $p\in\mathcal{P}^*(\R^n)$ and if there are constants $C>0$ and $p_\infty$ such that
\begin{equation}
|p(x)-p(y)|\le -\frac{C}{\log(|x-y|)},\hspace{1cm}x,y\in\R^n,\ ~|x-y|<1/2,
\end{equation}
and
\begin{equation}
|p(x)-p_\infty|\le \frac{C}{\log(e+|x|)},\hspace{1cm}x\in\R^n.
\end{equation}

The statement below establishes a reverse H\"older property under certain conditions on the exponents involved.
\begin{prop}(\cite[Lemma~2.7]{MP})\label{prop:revHolder}
Let $\p, r(\cdot)\in\mathcal{P}^{\log}(\R^n)$ such that $r(\cdot)\le\p$. Suppose
that $1/r(\cdot)=1/\p+1/s(\cdot)$. Then, for every ball $B\subset\R^{n}$
\begin{equation}\label{revHolder}
\|\chi_B\|_{r(\cdot)}\simeq\|\chi_B\|_{\p}\|\chi_B\|_{s(\cdot)}.
\end{equation}
\end{prop}
\
\begin{rem}\label{rem:revHolder}
It is straightforward to see that inequality~\eqref{revHolder} can be also written as
\begin{equation*}
\|\chi_B\|_{\pri}\simeq\|\chi_B\|_{r'(\cdot)}\|\chi_B\|_{s(\cdot)}.
\end{equation*}
\end{rem}

 Given a weight $w$ and $p\in\mathcal{P}(\R^n)$, we define the weighted variable Lebesgue space $L^{\p}(w)$ to be the set of all measurable functions $f$ such that $fw\in L^{\p}$, and we write 
\begin{equation*}
\|f\|_{L^{\p}(w)}=\|f\|_{\p,w}=\|fw\|_{\p}.
\end{equation*} 

Thus, we say that an operator $T$ is bounded on $L^{\p}(w)$ if 
\begin{equation*}
\|Tf\,w\|_{\p}\le C\|fw\|_{\p},
\end{equation*}
for all $f\in L^{\p}(w)$.

Following~\cite{Cabral} we define classes of variable weights associated with a function of critical radius $\rho$. Given a critical radius function $\rho$ and $p\in\mathcal{P}(\R^n)$ we  introduce the $A_{\p}^{\rho,\loc}$ the class of weights as those  $w$ for which there exists a constant $C>0$ such that the inequality
\begin{equation}\label{Ap_local}
\|w\chi_B\|_{\p}\|w^{-1}\chi_B\|_{p'(\cdot)}\le C|B|,
\end{equation}
holds for every ball $B\in\mathcal{B}_\rho$.

On the other hand, we will say that a weight $w$ belongs to the class $A_{\p}^{\rho}$ if there exist constants $\theta\ge 0$ and $C>0$ such that the inequality
\begin{equation}\label{Ap_rho}
\|w\chi_B\|_{\p}\|w^{-1}\chi_B\|_{p'(\cdot)}\le C|B|\left(1+\frac{r}{\rho(x)}\right)^\theta,
\end{equation}
holds for all balls $B=B(x,r)\subset \R^{n}$.

The following result shows the connection between the weights $A_{\p}^{\rho}$ and the boundedness of the operators $M_\rho^\theta$ in $L^{\p}(w)$.

\begin{thm}(\cite[Theorem~5]{Cabral})\label{teo:Mthetavar}
Let $p\in\mathcal{P}^{\log}(\R^n)$ with $p^{-}>1$. Then,  a weight $w\in A_{\p}^\rho$ if and only if there exists  $\theta>0$ such that $M^\theta_{\rho}$ is bounded on $L^{\p}(w)$.
\end{thm}

\begin{rem}\label{obs:pesos_var}
It follows from the above definition that if $w\in A_{\p}^\rho$, then $w^{-1}\in A_{\pri}^\rho$.
\end{rem}

Our second pair of results yields off-diagonal inequalities between two different weighted variable Lebesgue spaces. In the case of the constant exponent, the first of these was first demonstrated in~\cite{Pola}, while for the variable case, the same was done in~\cite{CUW} and~\cite{CUFMP} with and without weights, respectively.

To establish such a result we first define an appropriate class of weights that generalizes to the weights $A_{p,q}^\rho$.

Given a  critical radius function  $\rho$ and $\p, \q\in\mathcal{P}(\R^n)$, we say that a weight $w$ belongs to class $A_{\p,\q}^\rho$ if there exists a positive constant $C$ such that 
\begin{equation}\label{Apq_rho}
\|w\chi_B\|_{\q}\|w^{-1}\chi_B\|_{p'(\cdot)}\le C\|\chi_B\|_{\q}\|\chi_B\|_{p'(\cdot)}\left(1+\frac{r}{\rho(x)}\right)^\theta,
\end{equation}
for some $\theta>0$ and any ball $B=B(x,r)$ in $\R^n$. When $\p\in\mathcal{P}^{\log}(\R^n)$ and $\p=\q$ this is the $A_{\p}^\rho$ class.

\begin{rem}
If for some constant $0<\gamma<1$, $\frac1{p(x)}-\frac1{q(x)}=\gamma$, the condition~\eqref{Apq_rho} is equivalent to the following
\begin{equation*}\label{Apq_rhoal}
\|w\chi_B\|_{\q}\|w^{-1}\chi_B\|_{p'(\cdot)}\le C|B|^{1-\gamma}\left(1+\frac{r}{\rho(x)}\right)^\theta.
\end{equation*}
\end{rem}

This definition of class $A_{\p,\q}$ was adopted in~\cite{CUW} in the classical case, i.e. $V=0$.

We also say that $w\in A_{\p,\infty}^\rho$ if there exists a positive constant $C$ such
that
\begin{equation}\label{Apinfty}
\|w\chi_B\|_{\infty}\|w^{-1}\chi_B\|_{p'(\cdot)}\le C\|\chi_B\|_{p'(\cdot)}\left(1+\frac{r}{\rho(x)}\right)^\theta,
\end{equation}
for some $\theta>0$ and any ball $B=B(x,r)$ in $\R^n$.

We are now in a position to establish another main outcome of this section.

\begin{thm}\label{teo4}
Let $p_0$ and $q_0$ be such that $1<p_0<q_0<\infty$ and suppose that for all $w_0\in A_{p_0,q_0}^\rho$ it is verified that
\begin{equation}\label{hipoTeo4}
\|fw_0\|_{q_0}\le C\|gw_0\|_{p_0},\hspace{1cm} (f,g)\in\mathcal{F}.
\end{equation}
Be also  $\p,\q\in\mathcal{P}(\R^{n})$ such that
\begin{equation*}
\frac1{p(x)}-\frac1{q(x)}=\frac1{p_0}-\frac1{q_0}=\frac1{\sigma'}.
\end{equation*}

Then, if $w\in A_{\p,\q}^{\rho}$ and $\q\in\mathcal{P}^{\log}(\R^n)$ with $q^{-}>\sigma$, it follows that
\begin{equation*}
\|fw\|_{\q}\le C\|gw\|_{\p},\hspace{1cm} (f,g)\in\mathcal{F}.
\end{equation*}
\end{thm}
\

As a direct consequence of Theorems ~\ref{teo1} and~\ref{teo4} we obtain the following extrapolation result.

\begin{thm}\label{teo5}
Let $1<\beta<\infty$, $0\leq\delta<1$ and $\frac{1}{s}=\frac{1}{\beta}-\frac{\delta}{n}$ be such that the following inequality,
\begin{equation*}
\vertiii{f}_{\mathscr{L}_\delta(w)}\leq C \|g\,w\|_s,\ \ \ \ (f,g)\in\mathcal{F},
\end{equation*}
is verified for all $w \in A^{\rho}_{s,\infty}$ and some constant $C>0$.
Be also  $\p,\q\in\mathcal{P}(\R^{n})$ such that
\begin{equation*}
\frac1{p(x)}-\frac1{q(x)}=\frac1{\beta}.
\end{equation*}
	
Then, if $w\in A_{\p,\q}^{\rho}$ and $\q\in\mathcal{P}^{\log}(\R^n)$ with $q^{-}>\beta'$, it follows that
\begin{equation*}
\|fw\|_{\q}\le C\|gw\|_{\p},\hspace{1cm} (f,g)\in\mathcal{F}.
\end{equation*}
\end{thm}
\

To establish our next result, we will consider a variable version of the space $\mathscr{L}_\delta(w)$.

Let $w$ be a weight, $\gamma>1$ and $\p\in\mathcal{P}(\R^n)$ such that $\gamma\le p^{-}\le p(x)\le p^{+}<\frac{n\gamma}{(n-\gamma)_{+}}$ and let $\frac{\delta(x)}{n}=\frac1\gamma+\frac1{p(x)}$ (by $\beta_{+}$ we understand $\beta$ if $\beta>0$ and 0 is $\beta\le0$). The space $\mathscr{L}_{\delta(\cdot)}(w)$ is defined by the set of all functions $f$ on $L^{1}_{\loc}(\R^n)$ such that
\begin{equation}\label{oscilacionvar}
\frac{\|w\chi_B \|_\infty}{|B|^{\frac1{\gamma}}\|\chi_B\|_{\pri}}\int_B |f(y)-f_B|\,dy \leq C,\ \ \ \  \forall B \in \mathcal{B} _\rho.
\end{equation}  
and 
\begin{equation}\label{promediovar}
\frac{\|w\,\chi_{B(x,\rho(x))} \|_\infty}{|B(x,\rho(x))|^{\frac1{\gamma}}\|\chi_{B(x,\rho(x))}\|_{\pri}} \int_{B(x,\rho(x))}|f(y)|\,dy \leq C,\ \ \ \ \forall x\in\R^n.
\end{equation}
When $p(x)$ is equal to a constant $p$, this space coincides with the space $\mathscr{L}_{n/\gamma-n/p}(w)$ defined above.

We denote by $\vertiii{f}_{\mathscr{L}_{\delta(\cdot)}(w)}$ the smallest constant that satisfies \eqref{oscilacionvar} and \eqref{promediovar} simultaneously.

We are now in a position to state the last result of our interest.
\begin{thm}\label{teo6}
Let $\gamma>1$, $0\le\delta<\min\{n/\gamma,1\}$ and $s>1$ be such that $\frac1{s}=\frac1{\gamma}-\frac{\delta}{n}$. Let $\p\in\mathcal{P}^{\log}(\R^n)$ with $p^{-}>1$, and $0\le\tilde{\delta}(\cdot)<\delta$ such that $\frac1{\p}=\frac1{\gamma}-\frac{\tilde{\delta}(\cdot)}{n}$. Let us suppose that for every $w^{\eta}\in A^\rho_{\frac{s}{\eta},\infty}$, where $1<\eta<p^{-}$, it is verified that
\begin{equation*}
\vertiii{f}_{\mathscr{L}_{\delta}(w)}\leq C \|g\,w\|_s,\ \ \ \ (f,g)\in\mathcal{F},
\end{equation*}
Then,  for every $w^{\eta}\in A^\rho_{\frac{\p}{\eta},\infty}$ it holds that
\begin{equation*}
\vertiii{f}_{\mathscr{L}_{\tilde{\delta}(\cdot)}(w)}\leq C \|g\,w\|_{\p},\ \ \ \ (f,g)\in\mathcal{F}.
\end{equation*}
\end{thm}

\section{Auxiliary results}\label{Sec3}

Before proving the main results, we give several technical results necessary for proof. The first is a version of Rubio de Francia algorithm.
\begin{lem}\label{lem:alg RF cte}
Let $v \in A^{\rho}_m$ with $m \ge 1$. Then, for any $h\ge 0 $ belonging to $L^{m'} (v)$, there exists $H \in L^{m'}(v)$ such that
\begin{enumerate}
\item $h(x) \le H(x), \quad a.e.\  x \in \R^{n}$.
\item $\|H\|_{L^{m'}(v)} \leq C \|h\|_{L^{m'}(v)}$.
\item $Hv \in A^{\rho}_1$.
\end{enumerate}
\end{lem}

\begin{proof}
Following the so-called Rubio de Francia algorithm, it is sufficient to consider $H=\mathscr{R}h$ defined by
\begin{equation*}
\mathscr{R}h(x)=\sum_{k=0}^{\infty} \frac{[(M_\rho^{\theta})']^kh(x)}{2^k\|(M_\rho^{\theta})'\|^{k}_{L^{m'}(v)}},
\end{equation*}
where $(M_\rho^{\theta})'f=M_\rho^{\theta}(fv)/v$ and $[(M_\rho^{\theta})']^k$, for $k\ge1$ is $k$ times the composition of the operator $(M_\rho^{\theta})'$ and $[(M_\rho^{\theta})']^{0}$ denotes the identity operator.

Hence, the properties of $H$ follow easily from its definition and from bounding properties of $M_\rho^{\theta}$.

\end{proof}

We now define a localized version of the maximal sharp function as follows.

Given $f\in L^{1}_{\loc} (\R^{n})$ we define maximal sharp local function as
\begin{equation*}
f^{\sharp,\rho}_{\delta}(x)=\sup_{B\in\mathcal{B}_{\rho}} \frac{\chi_B(x)}{|B|^{1+\frac{\delta}{n}}} \int_B |f(y)-f_B|dy \quad + \sup_{B=B(z,\rho(z))} \frac{\chi_B(x)}{|B|^{1+\frac{\delta}{n}}} \int_B |f(y)|dy.
\end{equation*}

As expected, the space $\mathscr{L}_\delta(w)$ can be described by the operator $f^{\sharp,\rho}_\delta$ defined above.

\begin{lem}\label{lem: fsharp}
Let $0\leq\delta<1$ and $w$ be a weight. Then, there exist positive constants $C_1$ and $C_2$ such that
\begin{equation*}
C_1 \big\|w f^{\sharp,\rho}_\delta\big\|_\infty\leq \vertiii{f}_{\mathscr{L}_\delta(w)} \leq C_2 \big\|w f^{\sharp,\rho}_\delta\big\|_\infty.
\end{equation*}
\end{lem}
\begin{proof}
The proof of the lemma follows in a similar way to that given in Lemma 2 of \cite{BCH} for the case $\delta=0$.
\end{proof}

The following result is used in the proof of the Theorem~\ref{teo1}. It was proved in Corollary 5 of \cite{BCH}.

\begin{lem}\label{lem:desnormp}
Let $1<p<\infty$ and $w\in A^{\rho,\loc}_{s}$ for some $s\ge1$. If $g \in L^{1}_{\loc}\left( \R^{n}\right)$ then there exists a constant $C>0$ such that
\begin{equation*}
\|g\|_{L^p(w)} \leq C\|g_0^{\sharp,\rho}\|_{L^{p}(w)}.
\end{equation*}
\end{lem}

A key element in demonstrating our extrapolation results is the following generalization of a result given by Calder\'on and Scott in \cite[Proposition~4.6]{Calderon}.

\begin{prop}\label{prop:Calderon-Scott}
Let $0\leq\tilde{\delta}<\delta<1$, $1<p\le\frac{n}{\delta-\tilde{\delta}}$ y $\frac{1}{r}=\frac{1}{p} - \frac{\delta-\tilde{\delta}}{n}$. If $w \in A^{\rho,\loc}_{p,r}$ then there exists a positive constant $C$ such that
\begin{equation*}
\|f^{\sharp,\rho}_{\tilde{\delta}} w\|_r\leq C\|f^{\sharp,\rho}_{\delta}w\|_p.
\end{equation*}
\end{prop}

The proof of Proposition \ref{prop:Calderon-Scott} is a consequence of the following lemma which involves a localized fractional maximal function. 

Given a locally integrable function $f$ and $0\le\delta<1$,  $M^{\rho,\loc}_{\delta}$ is defined as follows,
\begin{equation*}
M^{\rho,\loc}_{\delta} f(x)=\sup_{B \in \mathcal{B}_{\rho}} \frac{\chi_{B}(x)}{|B|^{1-\frac{\delta}{n}}} \int_B|f(y)| dy.
\end{equation*}

\begin{lem}\label{lem:des CS}
Let $0\le\tilde{\delta}\le\delta<1$. Then,
\begin{equation}
f^{\sharp,\rho}_{\tilde{\delta}} (x) \leq 2 M^{\rho,\loc}_{\delta-\tilde{\delta}} \big(f^{\sharp,\rho}_{\delta}\big) (x).
\end{equation}
\end{lem}

\begin{proof}
Let $x\in\R^n$ and $B'\in\mathcal{B}_\rho$ such that $x\in B'$. Thus,
\begin{equation*}
\begin{split}
\frac{1}{|B'|^{1+\frac{\tilde{\delta}}{n}}} \int_{B'}|f(y)-f_{B'}|dy 
&=\frac{|B'|^{\frac{\delta-\tilde{\delta}}{n}}}{|B'|^{1+\frac{\delta-\tilde{\delta}}{n}}} \int_{B'} |f(y)-f_{B'}|\,dy\\
&\le|B'|^{\frac{\delta-\tilde{\delta}}{n}}f^{\sharp,\rho}_{\delta}(z), 
\end{split}
\end{equation*}
for any $z\in B'$. Integrating over $B'$ it follows that
\begin{equation*}
\begin{split}
|B'| \frac{1}{|B'|^{1+\frac{\tilde{\delta}}{n}}} \int_{B'}|f(y)-f_{B'}|\,dy 
&=\int_{B'}\left( \frac{1}{|B'|^{1+\frac{\tilde{\delta}}{n}}}\int_{B'}|f(y)-f_{B'}|\,dy \right) dz\\
&\leq|B'|^{\frac{\delta-\tilde{\delta}}{n}}\int_{B'}f^{\sharp,\rho}_{\delta}(z)\,dz,
\end{split}
\end{equation*}
from which it follows that
\begin{equation*}
\begin{split}
\frac{1}{|B'|^{1+\frac{\tilde{\delta}}{n}}} \int_{B'}|f(y)-f_{B'}|\,dy 
&\leq \frac{1}{|B'|^{1-\frac{\delta-\tilde{\delta}}{n}}}\int_{B'}f^{\sharp,\rho}_{\delta}(z)\,dz.\\
&\leq M^{\rho,\loc}_{\delta-\tilde{\delta}}\big(f^{\sharp,\rho}_{\delta}\big)(x).
\end{split}
\end{equation*}

Thus
\begin{equation}\label{sup1}
\sup_{B\in\mathcal{B}_{\rho}} \frac{\chi_{B}(x)}{|B|^{1+\frac{\tilde{\delta}}{n}}}\int_{B} |f(y)-f_{B}|\,dy\leq M^{\rho,\loc}_{\delta-\tilde{\delta}}\big(f^{\sharp,\rho}_{\delta}\big) (x).
\end{equation}

In a completely analogous manner, taking $B=B(z,\rho (z))$ with $z\in\R^n$ in place of $B'\in\mathcal{B}_\rho$, it is concluded that
\begin{equation}\label{sup2}
\sup_{B=B(z,\rho(z))} \frac{\chi_{B}(x)}{|B|^{1+\frac{\tilde{\delta}}{n}}} \int_{B} |f(y)|\,dy \leq M^{\rho,\loc}_{\delta-\tilde{\delta}}\big(f^{\sharp,\rho}_{\delta}\big)(x).
\end{equation}

From \eqref{sup1} and \eqref{sup2} it follows that
\begin{equation*}
f^{\sharp,\rho}_{\tilde{\delta}} (x) \leq 2 M^{\rho,\loc}_{\delta-\tilde{\delta}} \big(f^{\sharp,\rho}_{\delta}\big) (x).
\end{equation*}
\end{proof}

\begin{proof}[Proof of the Proposition~\ref{prop:Calderon-Scott}]
Suppose that $\| f_{\delta}^{\#,\rho} w\|_p < \infty$, that is, $f_{\delta}^{\#,\rho} \in L^{p}(w^p)$, otherwise, there is nothing to prove.

Let us first consider the case where $p<\frac{n}{\delta-\tilde{\delta}}$ and let $r$ be as in the hypothesis. By the previous lemma we have that
\begin{equation}\label{ftilde}
\begin{split}
\|f_{\tilde{\delta}}^{\#,\rho} w \|_r&=\left(\int_{\R^n}[f_{\tilde{\delta}}^{\#,\rho}w]^{r}\,dx \right)^{1/r}\\
&\le2\left(\int_{\R^d}\big[M_{\delta -\tilde{\delta}}^{\rho,\loc}\big(f_{\delta}^{\#,\rho}\big)\big]^{r}w^{r}\,dx  \right)^{1/r} \\
&=2\big\|M_{\delta-\tilde{\delta}}^{\rho,\loc}\big(f_{\delta}^{\#,\rho}\big)w\big\|_r.
\end{split}
\end{equation}

Considering  $\delta>\tilde{\delta}$ and $\alpha =\delta-\tilde{\delta}>0$ we have that  $\alpha<n$ and both $p$ and $r$ fulfill the hypotheses of the Theorem 3.1.13 in \cite{Cabral_T}, so that $M_{\delta-\tilde{\delta}}^{\rho,\loc}$ is bounded from $L^p(w^p)$ in $L^r (w^r)$.

Therefore, it follows from \eqref{ftilde} that
\begin{equation}
\|f_{\tilde{\delta}}^{\#,\rho} w \|_r \leq C\|f_{\delta}^{\#,\rho} w \|_p.
\end{equation}

On the other hand, if  $\delta=\tilde{\delta}$, it follows from the hypothesis that $r=p$. Also, $w \in A_{p,p}^{\rho}$ if and only if $w^p \in A_{p}^{\rho}$. Then, $w^{p} \in A_{p}^{\rho,\loc}$ and from Theorem 1 in~\cite{BHS} it follows that  $M_{\delta-\tilde{\delta}}^{\rho,\loc}$ is bounded on $L^{p}(w^{p})$, from which we obtain the result we wanted.

Let us now look at the case where  $p=\frac{n}{\delta-\tilde{\delta}}$, so that $r=\infty$. We want to see that if $w\in A_{\frac{n}{\delta-\tilde{\delta}},\infty}^{\rho,\loc}$, then there exists $C>0$ such that
\begin{equation*}
\|f_{\tilde{\delta}}^{\#,\rho} w \|_\infty \leq C\|f_{\delta}^{\#,\rho} w \|_{\frac{n}{\delta-\tilde{\delta}}}.
\end{equation*}

From Lemma \ref{lem:des CS} it follows that
\begin{equation} \label{rinf1}
\begin{split}
\|f_{\tilde{\delta}}^{\#,\rho} w\|_\infty
&\le2\|M_{\delta-\tilde{\delta}}^{\rho,\loc} (f_{\delta}^{\#,\rho}) w  \|_\infty.
\end{split}
\end{equation}

On the other hand, by Proposition 3.1.12 in~\cite{Cabral_T} it follows that if $w^{\frac{n}{\delta - \tilde{\delta}-n}} \in A_{1}^{\rho,\loc}$, then $M_{\delta-\tilde{\delta}}^{\rho,\loc}$ is bounded from $L^{\frac{n}{\delta- \tilde{\delta}}}(w^{\frac{n}{\delta -\tilde{\delta}}})$ on $L^{\infty}(w)$. It is enough to see that  $w^{\frac{n}{\delta-\tilde{\delta}-n}} \in A_{1}^{\rho,\loc}$. Being $p=\frac{n}{\delta-\tilde{\delta}}$, 
\begin{equation*}
\frac{n}{\delta-\tilde{\delta}-n}=\frac{1}{\frac{\delta-\tilde{\delta}}{n}-1}=\frac{1}{\frac{1}{p}-1}=
\frac{1}{-\frac{1}{p'}}=-p',
\end{equation*}
and since $w \in A_{p,\infty}^{\rho,\loc}$ it follows that $w^{-p'}\in A_{1}^{\rho,\loc}$, which is what we wanted to prove.

Therefore, by Proposition 3.1.12 in \cite{Cabral_T}, we obtain that
\begin{equation}\label{rinf2}
\|M_{\delta-\tilde{\delta}}^{\rho,\loc} (f_{\delta}^{\#,\rho})w \|_\infty \leq C \|f_{\delta}^{\#, \rho} w \|_{\frac{n}{\delta-\tilde{\delta}}}.
\end{equation}

Finally, from \eqref{rinf1} and \eqref{rinf2} we get that
\begin{equation*}
\|f_{\tilde{\delta}}^{\#,\rho}w \|_\infty \leq C \|f_{\delta}^{\#,\rho}w \|_{\frac{n}{\delta -\tilde{\delta}}}.
\end{equation*}

\end{proof}

We will also need some auxiliary results in the variable context.

The following results refer to properties of the weights  $A_{\p,\q}^{\rho}$. The proof of the first proposition is straightforward and we omit it.
\begin{prop}\label{propo:Apinf}
Let $\p,\q\in\mathcal{P}(\R^n)$. If $w\in A_{\p,\infty}^{\rho}$ then $w\in A_{\p,\q}^{\rho}$.
\end{prop}

\begin{prop}\label{propo:Apinfdec}
Let $\p,\q\in\mathcal{P}^{\log}(\R^n)$ such that $\q\le\p$. If $w\in A_{\q,\infty}^{\rho}$ then $w\in A_{\p,\infty}^{\rho}$.
\end{prop}
\begin{proof}
Let $s(\cdot)$ such that $\frac1{\q}=\frac1{\p}+\frac1{s(\cdot)}$.
Since $1/\pri=1/q'(\cdot)+1/s(\cdot)$, from H\"older's inequality, the hypothesis on $w$ and Remark~\ref{rem:revHolder}, we obtain that
\begin{equation*}
\begin{split}
\|w\chi_B\|_{\infty}\|w^{-1}\chi_B\|_{p'(\cdot)}&\le C\|w\chi_B\|_{\infty}\|w^{-1}\chi_B\|_{q'(\cdot)} \|\chi_B\|_{s(\cdot)} \\
&\le C\|\chi_B\|_{q'(\cdot)} \|\chi_B\|_{s(\cdot)} \left(1+\frac{r}{\rho(x)}\right)^\theta\\
&\le C\|\chi_B\|_{p'(\cdot)}\left(1+\frac{r}{\rho(x)}\right)^\theta.
\end{split}
\end{equation*}
\end{proof}

\begin{prop}\label{propo:Apqbis}
Let $\p, r(\cdot)\in\mathcal{P}^{\log}(\R^n)$ such that $\p\le r(\cdot)$, and $\eta$ satisfying $1<\eta<p^{-}$. If $w^\eta\in A_{\frac{\p}{\eta},\infty}^{\rho}$ then $w\in A^{\rho}_{r(\cdot),\infty}$.
\end{prop}
\begin{proof}
Since $\p\le r(\cdot)$ it follows from Proposition~\ref{propo:Apinfdec} that $w^\eta\in A_{\frac{r(\cdot)}{\eta},\infty}^{\rho}$. Taking $\beta=1/(\eta-1)$ it follows that $\frac{\eta}{r'(\cdot)}=\frac{1}{(r(\cdot)/\eta)'}+\frac1{\beta}$.

Then, by Proposition~\ref{prop:dilatacion}, the
generalized H\"older inequality and Proposition~\ref{prop:revHolder} we get
\begin{equation*}
\begin{split}
\left(\frac{\|w \chi_B \|_{\infty} \|w^{-1} \chi_B \|_{r'(\cdot)}}{\|\chi_B\|_{r'(\cdot)}}\right)^{\eta}
&=\frac{\|w^\eta \chi_B \|_{\infty} \|w^{-\eta} \chi_B \|_{r'(\cdot)/\eta}}{\|\chi_B\|_{r'(\cdot)/\eta}}\\
&\le C\frac{\|w^\eta \chi_B \|_{\infty} \|w^{-\eta} \chi_B \|_{(r(\cdot)/\eta)'}\|\chi_B\|_{\beta}}{\|\chi_B\|_{r'(\cdot)/\eta}}\\
&\le C\frac{\|w^\eta \chi_B \|_{\infty}\|w^{-\eta} \chi_B \|_{(r(\cdot)/\eta)'}}{\|\chi_B \|_{(r(\cdot)/\eta)'}}\\  
&\le C\left(1+\frac{r}{\rho(x)}\right)^{\theta}, 
\end{split}
\end{equation*}
where we have used $w^\eta\in A_{\frac{r(\cdot)}{\eta},\infty}^{\rho}$. We conclude that $w\in A_{r(\cdot),\infty}^{\rho}$.
\end{proof}

\begin{prop}\label{propo:Apqcaracbis}
Let $\p,\q\in\mathcal{P}^{\log}(\R^n)$ and $s>1$ such that $1<p^{-}\le\,p(x)<q(x)\le q^{+}<\infty$ and $\frac{1}{\p}-\frac{1}{\q}=\frac1{s'}$. Then, if $w\in A_{\p,\q}^{\rho}$   it is satisfied that $w^{-s}\in A^{\rho}_{\pri/s}$.
\end{prop}
\begin{proof}
Given a ball $B=B(x,r)\subset\R^n$, by Proposition~\ref{prop:dilatacion}, the inequality
\begin{equation*}
\|w^{-s}\chi_B\|_{\pri/s}\|w^{s}\chi_B\|_{(\pri/s)'}
\le C|B|\left(1+\frac{r}{\rho(x)}\right)^\theta.
\end{equation*}
is equivalent to
\begin{equation*}
\|w^{-1}\chi_B\|_{\pri}\|w\chi_B\|_{s(\pri/s)'}
\le C|B|^{1/s}\left(1+\frac{r}{\rho(x)}\right)^{\theta/s}.
\end{equation*}

We prove the last one. From Proposition~\ref{prop:revHolder},  since $\frac1{s}=\frac{1}{\q}+\frac{1}{\pri}$, we get that
\begin{equation*}
|B|^{1/s}=\|\chi_B\|_{s}\simeq\|\chi_B\|_{\q}\|\chi_B\|_{\pri}.
\end{equation*}
Then,  by the hypothesis on the weight $w$ and being $\q=s(\pri/s)'$ we obtain that 
\begin{equation*}
\begin{split}
\frac{\|w^{-1}\chi_B\|_{\pri}\|w\chi_B\|_{s(\pri/s)'}}{|B|^{1/s}}
&\le C\frac{\|w\chi_B\|_{\q}\|w^{-1}\chi_B\|_{\pri}}{\|\chi_B\|_{\q}\|\chi_B\|_{\pri}}\\
&\le C\left(1+\frac{r}{\rho(x)}\right)^{\vartheta}.
\end{split}
\end{equation*}
Then, taking $\theta=s\vartheta$, the proof is completed.

\end{proof}

\begin{prop}\label{propo:Apqcarac}
Let $\p,\q\in\mathcal{P}(\R^n)$ such that $1<p^{-}\le\,p(x)<q(x)\le q^{+}<\infty$ and we assume that it exists $\sigma>1$ such that $\frac{1}{\p}-\frac{1}{\q}=\frac1{\sigma'}$. Then, $w\in A_{\p,\q}^{\rho}$ if and only if $w^\sigma\in A^{\rho}_{\q/\sigma}$.
\end{prop}

\begin{rem}
Let us note that being $\frac{1}{\p}-\frac{1}{\q}=\frac1{\sigma'}$ we have $\frac{1}{p^{-}}-\frac{1}{q^{-}}=\frac1{\sigma'}$, that is,
\[
\frac{1}{p^{-}}-1=\frac{1}{q^{-}}-\frac1{\sigma},
\]
so the condition $p^{-}>1$ is equivalent to having $q^{-}>\sigma$.
\end{rem}

\begin{proof}
That weight  $w^\sigma$ belongs to $A_{\q /\sigma}^{\rho}$is equivalent to the existence of constants $\theta >0$ and $C>0$ such that
\begin{equation*}
\|w^{\sigma} \chi_B \|_{r(\cdot)}^{1/\sigma} \|w^{-\sigma} \chi_B \|_{r'(\cdot)}^{1/\sigma} \leq C |B|^{1/\sigma} \left(1+\frac{r}{\rho(x)}\right)^{\theta},
\end{equation*}
for any ball $B=B(x,r)$ in $\R^n$ and where $r(x)=q(x)/\sigma$.
	
Moreover, as $\frac{1}{p(x)}-\frac{1}{q(x)}=\frac{1}{\sigma '}$, by virtue of Proposition~\ref{prop:revHolder} $w$ belongs to  $A_{\p,\q}^{\rho}$if there are constants $C>0$ y $\theta'>0$ such that
\begin{equation*}
\begin{split}
\|w \chi_B \|_{\q} \|w^{-1} \chi_B \|_{p'(\cdot)} &\leq
C |B|^{1-\frac{1}{\sigma'}} \left(1+\frac{r}{\rho(x)}\right)^{\theta'}\\	 
&=  C |B|^{1/\sigma} \left(1+\frac{r}{\rho(x)}\right)^{\theta'}, 
\end{split}
\end{equation*}
for any ball $B=B(x,r)$ in $\R^n$. 
	
Since $1<q^{-}<q^{+}<\infty$, we can apply the Proposition~\ref{prop:dilatacion}  and obtain that for every ball $B$,
\begin{equation}\label{qs1}
\|w^{\sigma} \chi_B \|_{\q /\sigma}^{1/\sigma} =\|w \chi_B \|_{\q}.
\end{equation}
	
Now let us look at the relationship that exists between $r'(\cdot)$ and $p'(\cdot)$.
On the one hand, we have that 
\begin{equation}\label{rprima}
r'(x)=\frac{\frac{q(x)}{\sigma}}{\frac{q(x)}{\sigma}-1}=\frac{q(x)}{q(x)-\sigma}.
\end{equation}
On the other hand, 
\begin{equation*}
\frac{1}{p'(x)}=1-\frac{1}{p(x)}=1-\frac{1}{\sigma'}-\frac{1}{q(x)}=\frac{1}{\sigma}-\frac{1}{q(x)},
\end{equation*}
from which it follows that
\begin{equation}\label{rprima2}
p'(x)=\sigma \frac{q(x)}{q(x)-\sigma}
\end{equation}
	
Then, from \eqref{rprima} and \eqref{rprima2}, it follows that $r'(x)=\frac{1}{\sigma} p'(x)$. 
Then, again, by the Proposition~\ref{prop:dilatacion} and given that $(p')^{-}\ge1$ and $p^{-}>1$ we get that
\begin{equation}\label{qs2}
\|w^{-\sigma} \chi_B \|_{r'(\cdot)}^{1/\sigma} = \|w^{-1} \chi_B \|_{p'(\cdot)}.
\end{equation}

Therefore, from \eqref{qs1} and \eqref{qs2} it follows that $w\in A_{\p,\q}^{\rho}$ if, and only if, $w^{\sigma} \in A_{\q/\sigma}^{\rho}$.
	
\end{proof}

On the other hand, to construct the weight $w_0$ necessary to be able to apply the hypothesis of Theorem~\ref{teo4}, we will use item~\eqref{propie2} of Proposition~\ref{prop: prop_w} and therefore, to find the weights in $A_1^\rho$, we will consider the following variable version of Rubio de Francia's extrapolation algorithm.

\begin{prop}\label{prop:algorit_var}
Let $r(\cdot)\in\mathcal{P}(\R^n)$ and suppose that $v$ is a weight such that $M_\rho^\theta$ is bounded on  $L^{r(\cdot)}(v)$, for some $\theta>0$. For a non-negative function $h$ such that $h\in L^{r(\cdot)}(v)$ we define
\[
\mathcal{R}h(x)=\sum_{k=0}^{\infty}\frac{(M_\rho^{\theta})^kh(x)}{2^k\|M_\rho^{\theta}\|^k_{L^{r(\cdot)}(v)}}.
\]
	
Then,
\begin{enumerate}
\item\label{item1} $h(x)\le\mathcal{R}h(x)$,\ \ \ \ $a.e.\ x\in\R^n$.
\item\label{item2} $\|\mathcal{R}h\|_{L^{r(\cdot)}(v)}\le2\|h\|_{L^{r(\cdot)}(v)}$.
\item\label{item3} $\mathcal{R}h\in A_1^\rho$.
\end{enumerate}
\end{prop}
\begin{proof}
The proof is essentially the same as in the constant exponent case. Property $\eqref{item1}$ for $\mathcal{R}h$ is immediate, property $\eqref{item2}$ is deduced from the assumption that $M_\rho^{\theta}$ is bounded in $L^{r(\cdot)}(v)$, and finally, the property $\eqref{item3}$ follows from the fact that $M_\rho^{\theta}$ is sublinear and $h$ is non-negative.
	
\end{proof}

As a consequence of the previous result we will prove the following corollary which will be used in the proof of the Theorem~\ref{teo4}.

\begin{cor}\label{cor:algoritmo}
Let $\p\in\mathcal{P}(\R^n)$, $\q\in\mathcal{P}^{\log}(\R^n)$ and $\sigma$ such that $1<\sigma<q^{-}$, $1<p^{-}\le p(x)<q(x)\le q^{+}<\infty$  and $\frac1{\p}-\frac1{\q}=\frac1{\sigma'}$. Also let $h_1$ and $h_2$  be non-negative functions such that $h_1\in L^{\q}(w)$ and $h_2\in L^{(\q/\sigma)'}(\R^n)$.
	
Then, if $w\in A_{\p,\q}^{\rho}$ there exist $H_1\in L^{\q}(w)$ and $H_2\in L^{(\q/\sigma)'}(\R^n)$ those that verify,
\begin{multicols}{2}
\begin{enumerate}
\item\label{coritem1} $h_1(x)\le H_1(x)$,\ \ \ \ $a.e.\ x\in\R^n$.
\item\label{coritem2} $\|H_1\|_{L^{q(\cdot)}(w)}\le2\|h_1\|_{L^{q(\cdot)}(w)}$.
\item\label{coritem3} $H_1^{\,\sigma}\in A_1^\rho$.
\end{enumerate}
\begin{enumerate}[(1')]
\item\label{corpitem1} $h_2(x)\le H_2(x)$,\ \ \ \ $a.e.\ x\in\R^n$.
\item\label{corpitem2} $\|H_2\|_{(\q/\sigma)'}\le2\|h_2\|_{(\q/\sigma)'}$.
\item\label{corpitem3} $H_2\,w^\sigma\in A_1^\rho$.
\end{enumerate}
\end{multicols}
\end{cor}
\begin{proof}
Let us begin by proving the existence and properties of $H_1$. To do so, let us consider
\[
H_1=(\mathcal{R}\,h_1^\sigma)^{1/\sigma}.
\]
	
Since $w\in A_{\p,\q}^{\rho}$, in view of the Proposition~\ref{propo:Apqcarac}, we have that  $w^\sigma\in A^\rho_{\q/\sigma}$. Hence, by virtue of Theorem~\ref{teo:Mthetavar} and given that $q\in\mathcal{P}^{\log}(\R^n)$ with $q^{-}>\sigma$, one has that there exists $\theta>0$ such that $M_\rho^\theta$ is bounded in $L^{\q/\sigma}(w^\sigma)$.

Let us also note that since $h_1\in L^{\q}(w)$ it follows by Proposition~\ref{prop:dilatacion} that $h_1^{\sigma}\in L^{\q/\sigma}(w^\sigma)$.

Therefore, applying the Proposition~\ref{prop:algorit_var}, with $h=h_1^\sigma$, $r(\cdot)=\q/\sigma$ and $v=w^\sigma$, we have that
\begin{enumerate}[i)]
\item\label{coritemi} $h_1^\sigma\le\mathcal{R}\,h_1^\sigma$,\ \ \ \ $a.e$.
\item\label{coritemii} $\|\mathcal{R}\,h_1^\sigma\|_{L^{\q/\sigma}(w^\sigma)}\le2\|h_1^\sigma\|_{L^{\q/\sigma}(w^\sigma)}$.
\item\label{coritemiii} $\mathcal{R}\,h_1^\sigma\in A_1^\rho$.
\end{enumerate}
	
It is clear that~\textit{\ref{coritem1}.} and~\textit{\ref{coritem3}.} follow directly from ~\textit{\ref{coritemi})} and \textit{\ref{coritemiii})} respectively. On the other hand, ~\textit{\ref{coritem2}.} is a consequence of ~\textit{\ref{coritemii})}  and the Proposition~\ref{prop:dilatacion}, as
\begin{equation*}
\begin{split}
\|H_1w\|_{\q}&=\|(\mathcal{R}\,h_1^\sigma)^{1/\sigma}(w^\sigma)^{1/\sigma}\|_{\q}=\|\mathcal{R}h_1^\sigma\,w^\sigma\|^{1/\sigma}_{\q/\sigma}\\
&\le2^{1/\sigma}\|h_1^\sigma\, w^\sigma\|^{1/\sigma}_{\q/\sigma}
=2^{1/\sigma}\|h_1\,w\|_{\q}\le2\|h_1\,w\|_{\q}.
\end{split}
\end{equation*}
	
Let us now look at the existence and properties of $H_2$. We will take
\[
H_2=\mathcal{R}(h_2\,w^\sigma)w^{-\sigma}.
\]
	
Since $w\in A_{\p,\q}^{\rho}$, from the Proposition~\ref{propo:Apqcarac} and the Observation~\ref{obs:pesos_var}, we have that $w^{-\sigma}\in A^\rho_{(\q/\sigma)'}$. Thus, by Theorem~\ref{teo:Mthetavar} and given that $q\in\mathcal{P}^{\log}(\R^{n})$ and $(q^{+}/\sigma)'=[(\q/\sigma)']^{-}>1$ , one has that there exists $\theta'>0$ such that $M_\rho^{\theta'}$ is bounded in $L^{(\q/\sigma)'}(w^{-\sigma})$.
	
Let us also note that since  $h_2\in L^{(\q/\sigma)'}(\R^{n})$ it follows immediately that $h_2w^\sigma\in L^{(\q/\sigma)'}(w^{-\sigma})$.
	
Therefore, by applying Proposition~\ref{prop:algorit_var}, with $h=h_2w^\sigma$, $r(\cdot)=(\q/\sigma)'$ and $v=w^{-\sigma}$, we get
\begin{enumerate}[a)]
\item\label{coritema} $h_2w^\sigma\le\mathcal{R}(h_2w^\sigma)$,\ \ \ \ $a.e$.
\item\label{coritemb}  $\|\mathcal{R}(h_2w^\sigma)\|_{L^{(\q/\sigma)'}(w^{-\sigma})}\le2\|h_2w^\sigma\|_{L^{(\q/\sigma)'}(w^{-\sigma})}$.
\item\label{coritemc}  $\mathcal{R}(h_2w^\sigma)\in A_1^\rho$.
\end{enumerate}
	
Clearly, ~\textit{\ref{corpitem1}'.} and~\textit{\ref{corpitem3}'.} are a consequence of ~\textit{\ref{coritema})}  and~\textit{\ref{coritemc})}, respectively. Finally,~\textit{\ref{corpitem2}.} is obtained from~\textit{\ref{coritemb})} since
\begin{equation*}
\begin{split}
\|H_2\|_{(\q/\sigma)'}&=\|(\mathcal{R}h_2\,w^\sigma)w^{-\sigma}\|_{(\q/\sigma)'}\le2\|(h_2w^\sigma)w^{-\sigma}\|_{(\q/\sigma)'}=2\|h_2\|_{(\q/\sigma)'}.
\end{split}
\end{equation*}
	
\end{proof}

In the following corollary, given an exponent function $r(\cdot)$ and $\eta>1$, we denote by $\overline{r}=\overline{r}(\cdot)=\frac{r(\cdot)}{\eta}$.

\begin{cor}\label{cor:algoritmobis}
Let $\p,\q\in\mathcal{P}^{\log}(\R^n)$ and $s>1$ such that  $1<p^{-}\le p(x)<q(x)\le q^{+}<\infty$  and $\frac1{\p}-\frac1{\q}=\frac1{s}$. Let $1<\eta<p^{-}$ and non-negative function $\tilde{h}\in L^{\overline{q}/\overline{s}'}(w^{-\eta\overline{p}'\overline{s}'/\overline{q}})$.
	
Then, if $w^{\eta}\in A_{\p/\eta,\infty}^{\rho}$, there exists $\tilde{H}\in L^{\overline{q}/\overline{s}'}(w^{-\eta\overline{p}'\overline{s}'/\overline{q}})$ such that,
\begin{enumerate}
\item\label{corritem1} $\tilde{h}(x)\le\tilde{H}(x)$,\ \ \ \ $a.e.\ x\in\R^n$.
\item\label{corritem2} $\|\tilde{H}w^{-\eta\overline{p}'\overline{s}'/\overline{q}}\|_{\overline{q}/\overline{s}'}
\le2\|\tilde{h}w^{-\eta\overline{p}'\overline{s}'/\overline{q}}\|_{\overline{q}/\overline{s}'}$.
\item\label{corritem3} $\tilde{H}w^{-\eta\overline{p}'}\in A_1^\rho$.
\end{enumerate}
\end{cor}
\begin{proof}
Since $w^{\eta}\in A_{\overline{p},\infty}^{\rho}$, it follows from Proposition~\ref{propo:Apinf} that $w^{\eta}\in A_{\overline{p},\overline{q}}^{\rho}$. Moreover, given that $\frac1{\overline{p}(\cdot)}-\frac1{\overline{q}(\cdot)}=\frac1{\overline{s}}$, it holds by Proposition~\ref{propo:Apqcaracbis} that $w^{-\eta\overline{s}'}\in A_{\overline{p}'/\overline{s}'}^{\rho}$.
 	
Thus, the weight $w^{\eta\overline{s}'}\in A_{(\overline{p}'/\overline{s}')'}^{\rho}=A_{\overline{q}/\overline{s}'}^{\rho}$ and since that $p\in\mathcal{P}^{\log}(\R^n)$ with $p^{-}>1$, it follows that $(\overline{q}/\overline{s}')^{-}>1$, and by Theorem~\ref{teo:Mthetavar}  one has that there exists $\theta>0$ such that $M_\rho^{\theta}$ is bounded on $L^{\overline{q}/\overline{s}'}(w^{\eta\overline{s}'})$.

In addition, considering that $\tilde{h}w^{-\eta\overline{p}'}\in L^{\overline{q}/\overline{s}'}(w^{\eta\overline{s}'})$, we can apply the Proposition~\ref{prop:algorit_var} con $h=\tilde{h}w^{-\eta\overline{p}'}$, $r(\cdot)=\overline{q}/\overline{s}'$ and $v=w^{\eta\overline{s}'}$. We then have
\begin{enumerate}[i)]
\item\label{corritema} $\tilde{h}w^{-\eta\overline{p}'}\le\mathcal{R}(\tilde{h}w^{-\eta\overline{p}'})$,\ \ \ \ $a.e$.
\item\label{corritemb}  $\|\mathcal{R}(\tilde{h}w^{-\eta\overline{p}'})w^{\eta\overline{s}'}\|_{\overline{q}/\overline{s}'}
\le2\|\tilde{h}w^{-\eta\overline{p}'}w^{\eta\overline{s}'}\|_{\overline{q}/\overline{s}'}$.
\item\label{corritemc}  $\mathcal{R}(\tilde{h}w^{-\eta\overline{p}'})\in A_1^\rho$.
\end{enumerate}	

Thus, if we consider $\tilde{H}=\mathcal{R}(\tilde{h}w^{-\eta\overline{p}'})w^{\eta\overline{p}'}$ it is clear that~\textit{\ref{corritem1}.} and~\textit{\ref{corritem3}.} follow directly from ~\textit{\ref{corritema})} and \textit{\ref{corritemc})} respectively. Finally, ~\textit{\ref{corritem2}.} is a consequence of ~\textit{\ref{corritemb})} since
\begin{equation*}
\begin{split}
\|\tilde{H}w^{-\eta\overline{p}'\overline{s}'/\overline{q}}\|_{\overline{q}/\overline{s}'}&
=\|\mathcal{R}(\tilde{h}w^{-\eta\overline{p}'})w^{\eta\overline{p}'}w^{-\eta\overline{p}'\overline{s}'/\overline{q}}\|_{\overline{q}/\overline{s}'}\\
&=\|\mathcal{R}(\tilde{h}w^{-\eta\overline{p}'})w^{\eta\overline{s}'}\|_{\overline{q}/\overline{s}'}\\
&\le2\|\tilde{h}w^{-\eta\overline{p}'}w^{\eta\overline{s}'}\|_{\overline{q}/\overline{s}'}\\
&=2\|\tilde{h}w^{-\eta\overline{p}'\overline{s}'/\overline{q}}\|_{\overline{q}/\overline{s}'}.
\end{split}
\end{equation*}

\end{proof}

Similarly to the constant case, the space $\mathscr{L}_{\delta(\cdot)}(w)$ can be described by the opera-tor $f^{\sharp,\rho}_{\delta(\cdot),s(\cdot),\gamma}$ defined for $1<\gamma<\infty$ and $s(\cdot)\in\mathcal{P}(\R^n)$ such that $\delta(x)/n=1/\gamma-1/s(x)$ by
\begin{equation*}
\begin{split}
f^{\sharp,\rho}_{\delta(\cdot),s(\cdot),\gamma}(x)=\sup_{B\in\mathcal{B}_{\rho}} \frac{\chi_B(x)}{|B|^{1/\gamma}\|\chi_B\|_{s'(\cdot)}} &\int_B |f(y)-f_B|dy\\
&\ + \sup_{B=B(z,\rho(z))} \frac{\chi_B(x)}{|B|^{1/\gamma}\|\chi_B\|_{s'(\cdot)}} \int_B |f(y)|dy.
\end{split}
\end{equation*}

For convenience, in the rest of this paper,  we always assume that $f^{\sharp,\rho}_{\delta(\cdot),s(\cdot)}$ denotes $f^{\sharp,\rho}_{\delta(\cdot),s(\cdot),\gamma}$. The proof of the following lemma is analogous to the constant case.
	
\begin{lem}\label{lem: fsharpvar}
Let $0\leq\delta(\cdot)<1$ and $w$ be a weight. Then, there exist positive constants $C_1$ and $C_2$ such that
\begin{equation*}
C_1 \big\|w f^{\sharp,\rho}_{\delta(\cdot),s(\cdot)}\big\|_\infty\leq \vertiii{f}_{\mathscr{L}_{\delta(\cdot)}(w)} \leq C_2 \big\|w f^{\sharp,\rho}_{\delta(\cdot),s(\cdot)}\big\|_\infty.
\end{equation*}
\end{lem}

The following proposition is a key estimate in order to prove the Theorem~\ref{teo6}. A version in the classical case ($V=0$) can be found in~\cite{PRR}.
\begin{prop}\label{prop:CS-var}
Let $\p,s(\cdot)\in\mathcal{P}^{\log}(\R^n)$ and $\gamma>1$ such that  $\frac{\delta(\cdot)}{n}=\frac1{\gamma}-\frac1{s(\cdot)}$ and $\frac{\tilde{\delta}(\cdot)}{n}=\frac1{\gamma}-\frac1{\p}$ with $0\le\tilde{\delta}(\cdot)<\delta(\cdot)<1$. Then
\begin{equation*}
\big\|w f^{\sharp,\rho}_{\tilde{\delta}(\cdot),\,\p}\big\|_\infty\le C\big\|w f^{\sharp,\rho}_{\delta(\cdot),\,s(\cdot)}\big\|_{\frac{n}{\delta(\cdot)-\tilde{\delta}(\cdot)}},
\end{equation*}
for every $w\in A_{n/(\delta(\cdot)-\tilde{\delta}(\cdot)),\infty}^{\rho}$.
\end{prop}
\begin{proof}
Given $B=B(x_B,r_B)\in\mathcal{B}_\rho$, from definition of $f^{\sharp,\rho}_{\delta(\cdot),s(\cdot)}$, for almost everywhere $z\in B$, we have that
\begin{equation*}
\frac{1}{|B|^{1/\gamma}\|\chi_B\|_{s'(\cdot)}}\int_B |f(y)-f_B|dy\le f^{\sharp,\rho}_{\delta(\cdot),s(\cdot)}(z).
\end{equation*}
By integrating over $B$ and applying H\"older's inequality we obtain that
\begin{equation*}
\begin{split}
\frac{|B|}{|B|^{1/\gamma}\|\chi_B\|_{s'(\cdot)}}&\int_B |f(y)-f_B|dy\\
&\le \int_{B}f^{\sharp,\rho}_{\delta(\cdot),s(\cdot)}(z)\,dz\\
&\le C\big\|f^{\sharp,\rho}_{\delta(\cdot),s(\cdot)}\,w\big\|_{\frac{n}{\delta(\cdot)-\tilde{\delta}(\cdot)}}\big\|w^{-1}\chi_B\big\|_{\frac{n}{n-(\delta(\cdot)-\tilde{\delta}(\cdot))}},
\end{split}
\end{equation*}
from which it follows that
\begin{equation*}
\begin{split}
\frac{\|w\chi_B\|_{\infty}}{|B|^{1/\gamma}\|\chi_B\|_{p'(\cdot)}}&\int_B |f(y)-f_B|dy\\
&\le C\big\|f^{\sharp,\rho}_{\delta(\cdot),s(\cdot)}\,w\big\|_{\frac{n}{\delta(\cdot)-\tilde{\delta}(\cdot)}}\big\|w^{-1}\chi_B\big\|_{\frac{n}{n-(\delta(\cdot)-\tilde{\delta}(\cdot))}}
\frac{\|\chi_B\|_{s'(\cdot)}\|w\chi_B\|_{\infty}}{|B|\|\chi_B\|_{p'(\cdot)}}.
\end{split}
\end{equation*}

Considering the hypothesis on weight  $w$ and since $r_B\le\rho(x_B)$, we get that
\begin{equation*}
\begin{split}
\big\|w^{-1}\chi_B\big\|_{\frac{n}{n-(\delta(\cdot)-\tilde{\delta}(\cdot))}}&
\frac{\|\chi_B\|_{s'(\cdot)}\|w\chi_B\|_{\infty}}{|B|\|\chi_B\|_{p'(\cdot)}}\\
&\le C\big\|\chi_B\big\|_{\frac{n}{n-(\delta(\cdot)-\tilde{\delta}(\cdot))}}
\frac{\|\chi_B\|_{s'(\cdot)}}{|B|\|\chi_B\|_{p'(\cdot)}}\left(1+\frac{r_B}{\rho(x_B)}\right)^{\theta}\\
&\le C\big\|\chi_B\big\|_{\frac{n}{n-(\delta(\cdot)-\tilde{\delta}(\cdot))}}
\frac{\|\chi_B\|_{s'(\cdot)}}{|B|\|\chi_B\|_{p'(\cdot)}}.
\end{split}
\end{equation*}
Later, since $\frac{n-(\delta(\cdot)-\tilde{\delta}(\cdot))}{n}=\frac1{s(\cdot)}+\frac1{\pri}$, from Proposition~\ref{prop:revHolder} we have that
\begin{equation*}
\begin{split}
\big\|w^{-1}\chi_B\big\|_{\frac{n}{n-(\delta(\cdot)-\tilde{\delta}(\cdot))}}
\frac{\|\chi_B\|_{s'(\cdot)}\|w\chi_B\|_{\infty}}{|B|\|\chi_B\|_{p'(\cdot)}}&
\le C\|\chi_B\|_{s(\cdot)}\|\chi_B\|_{p'(\cdot)}
\frac{\|\chi_B\|_{s'(\cdot)}}{|B|\|\chi_B\|_{p'(\cdot)}}\\
&\le C\frac{\|\chi_B\|_{s(\cdot)}\|\chi_B\|_{s'(\cdot)}}{|B|}\\
&\le C.
\end{split}
\end{equation*}

Thus,
\begin{equation}\label{ec1}
\frac{\|w\chi_B\|_{\infty}}{|B|^{1/\gamma}\|\chi_B\|_{p'(\cdot)}}\int_B |f(y)-f_B|dy\le C\big\|f^{\sharp,\rho}_{\delta(\cdot),s(\cdot)
}\,w\big\|_{\frac{n}{\delta(\cdot)-\tilde{\delta}(\cdot)}}.
\end{equation}

On the other hand, considering $B'=B(t,\rho(t))$, we get
\begin{equation*}
\frac{1}{|B'|^{1/\gamma}\|\chi_{B'}\|_{s'(\cdot)}}\int_{B'} |f(y)|dy\le f^{\sharp,\rho}_{\delta(\cdot),s(\cdot)}(z),
\end{equation*}
for a.e. $z\in B'$. From here, proceeding analogously to the above, we get that
\begin{equation}\label{ec2}
\frac{\|w\chi_{B'}\|_{\infty}}{|B'|^{1/\gamma}\|\chi_{B'}\|_{p'(\cdot)}}\int_{B'} |f(y)|dy\le C\big\|f^{\sharp,\rho}_{\delta(\cdot),s(\cdot)}\,w\big\|_{\frac{n}{\delta(\cdot)-\tilde{\delta}(\cdot)}}.
\end{equation}
Finally, from~\eqref{ec1} and~\eqref{ec2}, the definition of $\vertiii{f}_{\mathscr{L}_{\tilde{\delta}(\cdot)}(w)}$ and Lemma~\ref{lem: fsharpvar}, the thesis follows.
\end{proof}

\section{Proof of main results}\label{Sec4}

In this section we give the proofs of our main theorems. We begin with the proof of Theorem~\ref{teo1}.

\begin{proof}[Proof of Theorem~\ref{teo1}]
Let $w \in A^{\rho}_{p,q}$ and $f \in L^{q}(w ^q)$. 
Without loss of generality, we can assume that  $\|g w\|_p< \infty$, since otherwise there is nothing to prove, and we can also assume that $\|g w\|_p>0$, since otherwise $g(x)=0$ at almost every point and then, from the hypothesis, we would also have $f(x)=0$ at almost every point.


Let us consider $h(x)= \left(\frac{|g(x)| w (x)^{p'}}{\|gw\|_p}\right)^{p-s}$, $\tilde{h} = h^{-s'/s}$ and $\frac{1}{r}=\frac{1}{p}-\frac{1}{s}$. Hence, $\tilde{h} \in L^{r/s'}(w^{-p'})$ and  	
\begin{equation*}
\int_{\R^n} \tilde{h}^{r/s'} w^{-p'}\,dx =1.
\end{equation*}

Since $w\in A^{\rho}_{p,q}$ and $\frac{1}{r}= \frac{1}{\beta}+ \frac{1}{q}-\frac{1}{s}=\frac{\delta}{n}+\frac{1}{q}$, it follows that $q\ge r$ and so $w\in A^{\rho}_{p,r}$. Therefore, by Proposition~\ref{prop: prop_w}, $w^{-p'}\in A^{\rho}_{1+\frac{p'}{r}}$ and applying Lemma~\ref{lem:alg RF cte}, with $m=1+\frac{p'}{r}$ and $v=w^{-p'}$, we know it exists $\tilde{H} \in L^{r/s'}(w^{-p'})$ such that $\tilde{H}\ge \tilde{h}$ and
\begin{equation}
\begin{split}\label{*}
\big\|\tilde{H}^{1/s'} w^{-p'/r}\big\|^r_r &= \int_{\R^{n}} \tilde{H}^{r/s'} w^{-p'}\,dx\\
&\leq C \int_{\R^n}\tilde{h}^{r/s'}w^{-p'}\,dx\\ 
&=C.
\end{split}	
\end{equation}

Hence,
\begin{equation*}
\begin{split}
\left(\int_{\R^n}|g|^p w^p dx\right)^{1/p} &= \left( \int_{\R^n} |gw^{p'}|^s h w^{-p'}\,dx  \right)^{1/s}\\
&=\left(\int_{\R^n} |g w^{p'}|^s(\tilde{h} ^{-1/s'})^sw^{-p'} dx\right)^{1/s}\\
&\ge\left(\int_{\R^n}|g w^{p'}|^s(\tilde{H} ^{-1/s'})^sw^{-p'} dx\right)^{1/s}\\
&=\left(\int_{\R^n}|g|^s(\tilde{H} ^{-1/s'}w^{p'/s'})^s dx\right)^{1/s}.
\end{split}
\end{equation*}

From item 3 of the Lemma~\ref{lem:alg RF cte}, $\tilde{H} w^{-p'} \in A^{\rho}_1$ and so, $\tilde{H}^{-1/s'} w^{p'/s'} \in A^{\rho}_{s,\infty}$. Thus, based on the hypothesis, the Lemma~\ref{lem: fsharp} and (\ref{*}) we obtain
\begin{equation}\nonumber
\begin{split}
\left(\int_{\R^n}|g|^p w^p\,dx\right)^{1/p}&\ge\left(\int_{\R^n}|g|^s\left(\tilde{H} ^{-1/s'}w^{p'/s'} \right)^s\,dx\right)^{1/s}\\
&\ge C\vertiii{f}_{\mathscr{L}_{\delta}(\tilde{H} ^{-1/s'}w^{p'/s'})}\\
& \ge C\|f^{\sharp,\rho}_{\delta}\tilde{H} ^{-1/s'}w^{p'/s'}\|_{\infty}\\
& \ge C\|f^{\sharp,\rho}_{\delta}\tilde{H} ^{-1/s'}w^{p'/s'}\|_{\infty}\|\tilde{H} ^{1/s'}w^{-p'/r}\|_r \\
& \ge C\|f^{\sharp,\rho}_{\delta} w\|_r.
\end{split}
\end{equation}

Moreover, considering that $w \in A^{\rho}_{p,q}$ and $r>p$, we obtain that $w \in A^{\rho}_{r,q}$. Then, by Proposition~\ref{prop:Calderon-Scott} and the fact that $A^{\rho}_{r,q}\subset A^{\rho,\loc}_{r,q}$, we have	
\begin{equation}
\left(\int_{\R^n}|g|^{p} w^{p}\,dx\right)^{1/p}\ge C\|f^{\sharp,\rho}_{\delta} w\|_r \ge C \|f^{\sharp,\rho}_{0} w\|_q
\end{equation} 
because $\frac{1}{q}= \frac{1}{r}-\frac{\delta}{n}$.

Finally, since $w\in A^{\rho}_{p,q}$, it follows that $w^{q} \in A^{\rho}_{1+\frac{q}{p'}} \subset A^{\rho, \loc}_{1+\frac{q}{p'}}$ and thus, by the Lemma~\ref{lem:desnormp}, we get that
\begin{equation*}
\|f^{\sharp,\rho}_0 w\|_q \ge C \|f w\|_q,
\end{equation*}
which concludes the proof.

\end{proof}

\begin{proof}[Proof of Theorem~\ref{teo2}]
Let  $w\in A^{\rho}_{p,\infty}$ and suppose as before that $\|gw\|_{p}<\infty$. Let us take $r$ such that  $\frac{1}{r}=\frac{1}{p}-\frac{1}{s}=\frac{\delta-\tilde{\delta}}{n}$. 

By items 3. and 1. of Proposition~\ref{prop: prop_w}, since $w\in A^\rho_{p,\infty}$, $w^{-p'}\in A^\rho_{1}$ and thus $w^{-p'}\in A^\rho_{1+p'/r}$.

In the same way as for the previous theorem, given $g$ we obtain a function $\tilde{H}\in L^{r/s'}(w^{-p'})$ such that
\begin{equation*}
\|\tilde{H}^{1/s'}w^{-p'/r}\|_{r}\leq C,
\end{equation*}
with $\tilde{H}^{-1/s'}w^{p'/s'}\in A^{\rho}_{s,\infty}$ and thereafter, from the hypothesis
\begin{equation*}
\left(\int_{\R^n}|g|^{p}w^{p}\,dx\right)^{1/p}\geq C\|f^{\sharp,\rho}_{\delta}w\|_{r}.
\end{equation*}

Now, since $w\in A^\rho_{p,\infty}$ and $r\ge p$, by item 6. of Proposition~\ref{prop: prop_w} we have $w\in A^{\rho}_{r,\infty}$.


Then, since $A^{\rho}_{r,\infty}\subset A^{\rho,\loc}_{r,\infty}$, applying Proposition~\ref{prop:Calderon-Scott} and Lemma~\ref{lem: fsharp} we get that
\begin{equation*}
\left(\int_{\R^n}|g|^{p}w^{p}\,dx\right)^{1/p}\geq C\|f^{\sharp,\rho}_{\delta}w\|_{r} \geq C
\|f^{\sharp,\rho}_{\tilde{\delta}}w\|_{\infty}\ge\vertiii{f}_{\mathscr{L}_{\tilde{\delta}}(w)}.
\end{equation*}
\end{proof}

\begin{proof}[Proof of Theorem~\ref{teo4}]
Let $f$ such that $\|f\|_{L^{\q}(w)}<\infty$ and suppose as before that $0<\|g\|_{L^{\p}(w)}<\infty$.
%
	
Let us take
\[
h_1=\frac{f}{\|fw\|_{\q}}+{w}^{-1}\left(\frac{g\,w}{\|gw\|_{\p}}\right)^{\frac{\p}{\q}}.
\]
	
Let us first see that $\|h_1w\|_{\q}\le C$. Indeed, from the definition of the norm, it follows that
\begin{equation*}
\begin{split}
\varrho_{\q}(h_1w)&=\int_{\R^n}|h_1(x)w(x)|^{q(x)}dx=\int_{\R^n}\left|\frac{f(x)w(x)}{\|fw\|_{\q}}+\left(\frac{g(x)w(x)}{\|gw\|_{\p}}\right)^{\frac{p(x)}{q(x)}}\right|^{q(x)}dx\\
&\le2^{q^{+}}\int_{\R^n}\left(\frac{|f(x)w(x)|}{\|fw\|_{\q}}\right)^{q(x)}dx+2^{q^{+}}\int_{\R^n}\left(\frac{|g(x)w(x)|}{\|gw\|_{\p}}\right)^{p(x)}dx\\
&\le2^{q^{+}+1}.
\end{split}
\end{equation*}
So, considering the Proposition~\ref{prop: modular} it follows what has been stated.
	
We will now use the Corollary~\ref{cor:algoritmo} and the operators $H_1$ and $H_2$ defined there for $h_1$ and the following $h_2$. Let $r_0=q_0/\sigma$. By Proposition~\ref{prop:dualidad}, there exists $h_2\in L^{(\q/\sigma)'}(\R^n)$, $h_2\ge0$, $\|h_2\|_{(\q/\sigma)'}\le1$, such that for any $\gamma>0$,
\begin{equation*}
\begin{split}
\|fw\|_{\q}^\sigma&=\|f^\sigma w^\sigma\|_{\q/\sigma}\le C\int_{\R^n}f^\sigma w^\sigma h_2\,dx
\le C\int_{\R^n}f^\sigma H_1^{-\gamma} H_1^{\gamma}w^\sigma H_2\,dx\\
&\le C\left(\int_{\R^n}f^{q_0} H_1^{-\gamma r_0}H_2 w^\sigma\,dx\right)^{1/r_0}
\left(\int_{\R^n}H_1^{\gamma r_0'}H_2 w^\sigma\,dx\right)^{1/r_0'}\\
&=CA^{1/r_0}B^{1/r_0'},
\end{split}
\end{equation*}
where we have consecutively considered the Proposition~\ref{prop:dilatacion}, the pointwise inequa-lity between $h_2$ and $H_2$ and the H\"older inequality with exponent $r_0$ and measure  $w^\sigma H_2\,dx$.
	
Let us now verify that the factor $B$ is uniformly bounded. On the one hand, considering the variable H\"older's inequality  (Proposition~\ref{prop:Holder}) with exponent $\q/\sigma$, we have that
\begin{equation*}
B\le C\|H_1^{\gamma r_0'}w^\sigma\|_{\q/\sigma}\|H_2\|_{(\q/\sigma)'}.
\end{equation*}

Thus, since $h_1\in L^{\q}(w)$ and $h_2\in L^{(\q/\sigma)'}(\R^n)$, taking $\gamma=\frac{\sigma}{r_0'}$ and considering again the Corollary~\ref{cor:algoritmo} and the Proposition~\ref{prop:dilatacion} we can see that  
\begin{equation*}
\|H_1^{\gamma r_0'}w^\sigma\|_{\q/\sigma}=\|H_1^{\sigma}w^\sigma\|_{\q/\sigma}=\|H_1w\|^\sigma_{\q}\le2^\sigma\|h_1w\|^\sigma_{\q}\le C,
\end{equation*}
and
\begin{equation*}
\|H_2\|_{(\q/\sigma)'}\le2\|h_2\|_{(\q/\sigma)'}\le C.
\end{equation*}

Let us now observe that from the choice of $\gamma$, we have that
\[
A=\int_{\R^n}f^{q_0}H_1^{-q_0/r_0'}H_2 w^\sigma\,dx.
\]
	
In order to use the hypothesis~\eqref{hipoTeo4}, we must jointly prove that $A$ is finite and that the weight $w_0=\big(H_1^{-q_0/r_0'}H_2 w^\sigma\big)^{1/q_0}$ belongs to $A_{p_0,q_0}^\rho$.
	
Let's start by seeing that $A$ is finite. Considering the definition of $h_1$, the properties of $H_1$ and $H_2$ and the variable H\"older's inequality it follows that
\begin{equation*}
\begin{split}
A&=\int_{\R^n}f^{q_0}H_1^{-q_0}H_1^{q_0}H_1^{-q_0/r_0'}H_2 w^\sigma\,dx\\
&\le\int_{\R^n}f^{q_0}h_1^{-q_0}H_1^{q_0}H_1^{-q_0/r_0'}H_2 w^\sigma\,dx\\
&\le\int_{\R^n}f^{q_0}\left(\frac{f}{\|fw\|_{\q}}\right)^{-q_0}H_1^{q_0}H_1^{-q_0/r_0'}H_2 w^\sigma\,dx\\
&=\|fw\|^{q_0}_{\q}\int_{\R^n}H_1^{q_0}H_1^{-q_0/r_0'}H_2 w^\sigma\,dx\\
&=\|fw\|^{q_0}_{\q}\int_{\R^n}H_1^{\sigma}H_2 w^\sigma\,dx\\
&\le C\|fw\|^{q_0}_{\q}\|H_1^\sigma w^\sigma\|_{\q/\sigma}\|H_2\|_{(\q/\sigma)'}<\infty.
\end{split}
\end{equation*}
	
Let us now see that  $w_0=\big(H_1^{-q_0/r_0'}H_2 w^\sigma\big)^{1/q_0}\in A_{p_0,q_0}^\rho$. For this purpose, considering the \'item~\textit{\ref{propie4}.} of the Proposition~\ref{prop: prop_w}, it is sufficient to prove that $w_0^{q_0}=H_1^{-(q_0-\sigma)}H_2 w^\sigma\in A_{s_0}^\rho$, where
\begin{equation*}
s_0=1+\frac{q_0}{p_0'}=\frac{q_0}{\sigma}.
\end{equation*}

In view again of Proposition~\ref{prop: prop_w}, but in this case \'item~\textit{\ref{propie2}.}, given that both $H_2 w^\sigma$ and		 $H_1^{\frac{q_0-\sigma}{s_0-1}}=H_1^\sigma$ belong to $A_1^\rho$, it follows immediately that,  
\begin{equation*}
w_0^{q_0}=\Big(H_1^{\frac{q_0-\sigma}{s_0-1}}\Big)^{1-s_0}H_2 w^\sigma\in A_{s_0}^\rho.
\end{equation*}
	
We are now in a position to apply~\eqref{hipoTeo4}. Considering then the hypothesis, the definition of $h_1$ and H\"older's inequality with respect to a certain exponent $r(\cdot)$, we have that
\begin{equation*}
\begin{split}
A^{1/r_0}=A^{\sigma/q_0}&=\left(\int_{\R^n}f^{q_0}H_1^{-q_0/r_0'}H_2 w^\sigma\,dx\right)^{\sigma/q_0}\\
&\le C\left(\int_{\R^n}g^{p_0}\big(H_1^{-q_0/r_0'}H_2 w^\sigma\big)^{p_0/q_0}dx\right)^{\sigma/p_0}\\
&\le C\|gw\|^\sigma_{\p}\left(\int_{\R^n}\bigg(h_1^{\frac{\q}{\p}}w^{\frac{\q}{\p}-1}\bigg)^{p_0}H_1^{-p_0/r_0'}H_2^{p_0/q_0} w^{\sigma p_0/q_0} dx\right)^{\sigma/p_0}\\
&\le C\|gw\|^\sigma_{\p}\left(\int_{\R^n}H_1^{p_0\big(\frac{\q}{\p}-\frac1{r_0'}\big)} w^{p_0\big(\frac{\sigma}{q_0}+\frac{\q}{\p}-1\big)}H_2^{p_0/q_0} dx\right)^{\sigma/p_0}\\
&\le C\|gw\|^\sigma_{\p}\bigg\|(H_1\,w)^{p_0\big(\frac{\q}{\p}-\frac1{r_0'}\big)}\bigg\|^{\sigma/p_0}_{r'(\cdot)}\ \bigg\|H_2^{p_0/q_0}\bigg\|^{\sigma/p_0}_{r(\cdot)}\\
&=C\|gw\|^\sigma_{\p}E^{\sigma/p_0}F^{\sigma/p_0},
\end{split}
\end{equation*}
where we have used the fact that $\frac{\sigma}{q_0}-1=-\frac1{r_0'}$.
	
Let us take $r(\cdot)=\frac{q_0}{p_0}(\frac{\q}{\sigma})'$ and see that both $E$ and $F$ are bounded. On the one hand, considering the Proposition~\ref{prop:dilatacion}, we get that  
\begin{equation*}
F=\Big\|H_2^{p_0/q_0}\Big\|_{r(\cdot)}=\Big\|H_2^{p_0/q_0}\Big\|_{\frac{q_0}{p_0}(\frac{\q}{\sigma})'}
=\big\|H_2\big\|^{p_0/q_0}_{(\q/\sigma)'}\le2^{p_0/q_0}\big\|h_2\big\|^{p_0/q_0}_{(\q/\sigma)'}\le2^{p_0/q_0}.
\end{equation*}
	
On the other hand, since
\begin{equation*}
p_0\bigg(\frac{\q}{\p}-\frac1{r_0'}\bigg)r'(\cdot)=\q,
\end{equation*}
we have that
\begin{equation*}
\begin{split}
\varrho_{r'(\cdot)}\bigg((H_1\,w)^{p_0\big(\frac{\q}{\p}-\frac1{r_0'}\big)}\bigg)
&=\int_{\R^n}(H_1w)^{p_0\big(\frac{\q}{\p}-\frac1{r_0'}\big)r'(\cdot)}\,dx\\
&=\int_{\R^n}(H_1w)^{\q}dx=\varrho_{\q}(H_1w).
\end{split}
\end{equation*}
	
Given that $\|H_1w\|_{\q}\le2\|h_1w\|_{\q}$, applying Proposition~\ref{prop: modular} it follows that $\varrho_{\q}(H_1w)$ is uniformly bounded, and therefore, applying Proposition~\ref{prop: modular} again, the uniform bounding of $E$ follows, which concludes the proof of the theorem.
	
\end{proof}

\begin{proof}[Proof of Theorem~\ref{teo6}]
Suppose as before that $0<\|g\|_{L^{\p}(w)}<\infty$.

Let $\frac1{\q}=\frac1{\p}-\frac1s=\frac{\ \delta-\tilde{\delta}(\cdot)}{n}$, $\overline{q}=\overline{q}(\cdot)=\frac{\q}{\eta}$, $\overline{p}=\overline{p}(\cdot)=\frac{\p}{\eta}$ and $\overline{s}=\frac{s}{\eta}$.\\

Let us consider the functions $h$ and $\tilde{h}$ defined by
\begin{equation*}
h=\left(\frac{|g|^\eta w^{\eta\overline{p}'}}{\|gw\|^\eta_{\p}}\right)^{\overline{p}-\overline{s}}
=\left(\frac{|g| w^{\overline{p}'}}{\|gw\|_{\p}}\right)^{\p-s}
\ \ \ \ \text{and}\ \ \ \ \tilde{h}=h^{-\overline{s}'/\overline{s}}.
\end{equation*}

From its definition, it follows that, $h\in L^{(\p/s)'}(w^{-\eta\overline{p}'/(\p/s)'})$ and \begin{equation*}
\|h\,w^{-\eta\overline{p}'/(\p/s)'}\|_{(\p/s)'}\le1.
\end{equation*} 
Indeed, since $\p=\Big(\frac{\p}{s}\Big)'(\p-s)=\overline{p}'(\p-\eta)$ we have that
\begin{equation*}
\begin{split}
\int_{\R^n}h(x)^{(p(x)/s)'}w(x)^{-\eta\overline{p}'}dx&=\int_{\R^n}\left(\frac{|g(x)| w(x)^{\overline{p}'}}{\|gw\|_{\p}}\right)^{p(x)}w(x)^{-\eta\overline{p}'}dx\\
&=\int_{\R^n}\frac{|g(x)w(x)|^{p(x)}}{\|gw\|^{p(x)}_{\p}}dx\le1.
\end{split}
\end{equation*}

Thus,
\begin{equation*}
\begin{split}
\int_{\R^n}\tilde{h}(x)^{\overline{q}/\overline{s}'}w(x)^{-\eta\overline{p}'}dx&=\int_{\R^n}h(x)^{-q(x)/s}w(x)^{-\eta\overline{p}'}dx\\
&=\int_{\R^n}h(x)^{(p(x)/s)'}w(x)^{-\eta\overline{p}'}dx\le1,
\end{split}
\end{equation*}
that is $\tilde{h}\in L^{\overline{q}/\overline{s}'}(w^{-\eta\overline{p}'\overline{s}'/\overline{q}})$ and $\|\tilde{h}w^{-\eta\overline{p}'\overline{s}'/\overline{q}}\|_{\overline{q}/\overline{s}'}\le1$.
\\

Moreover,
\begin{equation}\label{ecu}
\begin{split}
\int_{\R^n}|g(x)^\eta w(x)^{\eta\overline{p}'}|^{\overline{s}}\tilde{h}(x)^{-\overline{s}/\overline{s}'}w(x)^{-\eta\overline{p}'}dx
&=\int_{\R^n}|g(x)^\eta w(x)^{\eta\overline{p}'}|^{\overline{s}}h(x)w(x)^{-\eta\overline{p}'}dx\\
&=\int_{\R^n}\frac{(|g(x)| w(x)^{\overline{p}'})^{p(x)}}{\|gw\|^{p(x)-s}_{\p}}w(x)^{-\eta\overline{p}'}dx\\
&=\int_{\R^n}\frac{(|g(x)| w(x))^{p(x)}}{\|gw\|^{p(x)-s}_{\p}}dx\\
&\le\|gw\|_{\p}^s.
\end{split}
\end{equation}

Note that, since $\frac1{\q}=\frac1{\p}-\frac1{s}$, $w^\eta\in A_{\p/\eta,\infty}^\rho$ and $\tilde{h}\in L^{\overline{q}/\overline{s}'}(w^{-\eta\overline{p}'\overline{s}'/\overline{q}})$, by Corollary~\ref{cor:algoritmobis} there exists $\tilde{H}\in L^{\overline{q}/\overline{s}'}(w^{-\eta\overline{p}'\overline{s}'/\overline{q}})$ such that,
\begin{enumerate}
\item\label{corritem1b} $\tilde{h}(x)\le\tilde{H}(x)$,\ \ \ \ $a.e.\ x\in\R^n$.
\item\label{corritem2b} $\|\tilde{H}w^{-\eta\overline{p}'\overline{s}'/\overline{q}}\|_{\overline{q}/\overline{s}'}
\le2\|\tilde{h}w^{-\eta\overline{p}'\overline{s}'/\overline{q}}\|_{\overline{q}/\overline{s}'}$.
\item\label{corritem3b} $\tilde{H}w^{-\eta\overline{p}'}\in A_1^\rho$.
\end{enumerate}

Note that from item~\ref{corritem2b} and the fact that $\|\tilde{h}w^{-\eta\overline{p}'\overline{s}'/\overline{q}}\|_{\overline{q}/\overline{s}'}\le1$, it follows that $\|\tilde{H}w^{-\eta\overline{p}'\overline{s}'/\overline{q}}\|_{\overline{q}/\overline{s}'}\lesssim C$, which by Proposition~\ref{prop:dilatacion}, is equivalent to
\begin{equation}\label{ecub}
\|\tilde{H}^{1/(\eta\overline{s}')}w^{-\overline{p}'/\overline{q}}\|_{\q}\lesssim C.
\end{equation}

From~\eqref{ecu} and item~\ref{corritem1b} we obtain that
\begin{equation*}
\begin{split}
\|gw\|_{\p}&\ge\left(\int_{\R^n}|g(x)^\eta w(x)^{\eta\overline{p}'}|^{\overline{s}}\tilde{H}(x)^{-\overline{s}/\overline{s}'}w(x)^{-\eta\overline{p}'}dx\right)^{1/s}\\
&\ge\left(\int_{\R^n}|g(x)|^s\Big(\tilde{H}(x)^{-1/(\eta\overline{s}')}w(x)^{\overline{p}'/\overline{s}'}\Big)^{s}dx\right)^{1/s}\\
\end{split}
\end{equation*}

On the other hand, from item~\ref{corritem3b} we get that the weight $(\tilde{H}^{-1/\overline{s}'}w^{\eta\overline{p}'/\overline{s}'})^{-\overline{s}'}=\tilde{H}w^{-\eta\overline{p}'}\in A_1^\rho$, which is equivalent, by the item~\ref{propie3}. of Proposition~\ref{prop: prop_w}, to $\tilde{H}^{-1/\overline{s}'}w^{\eta\overline{p}'/\overline{s}'}\in A_{\overline{s},\infty}^\rho$, that is $(\tilde{H}^{-1/(\eta\overline{s}')}w^{\overline{p}'/\overline{s}'})^{\eta}\in A_{s/\eta,\infty}^\rho$. Thus, by the hypothesis, Lemma~\ref{lem: fsharp} and~\eqref{ecub} we obtain
\begin{equation*}
\begin{split}
\|gw\|_{\p}&\ge C\vertiii{f}_{\mathscr{L}_{\delta}(\tilde{H}^{-1/(\eta\overline{s}')}w^{\overline{p}'/\overline{s}'})}\\
&\ge C\|f_\delta^{\sharp,\rho}\tilde{H}^{-1/(\eta\overline{s}')}w^{\overline{p}'/\overline{s}'}\|_\infty\\
&\ge C\|f_\delta^{\sharp,\rho}\tilde{H}^{-1/(\eta\overline{s}')}w^{\overline{p}'/\overline{s}'}\|_\infty\|\tilde{H}^{1/(\eta\overline{s}')}w^{-\overline{p}'/\overline{q}}\|_{\q}\\
&\ge C\|f_\delta^{\sharp,\rho}w^{\overline{p}'(1/\overline{s}'-1/\overline{q})}\|_{\q}\\
&=C\|f_\delta^{\sharp,\rho}w\|_{\q}
\end{split}
\end{equation*}

From the hypothesis on $w$ and Proposition~\ref{propo:Apqbis} we conclude that $w\in A_{\q,\infty}^\rho$ with $\q=\frac{n}{\delta-\tilde{\delta}(\cdot)}$. Thus, from Proposition~\ref{prop:CS-var} and Lemma~\ref{lem: fsharpvar} we get that
\begin{equation*}
\begin{split}
\|gw\|_{\p}&\ge C\|f_\delta^{\sharp,\rho}w\|_{\q}\\
&=C\|f_\delta^{\sharp,\rho}w\|_{\frac{n}{\delta-\tilde{\delta}(\cdot)}}\\
&\ge C\|f_{\tilde{\delta}(\cdot),\,\p}^{\sharp,\rho}w\|_{\infty}\\
&\ge C\vertiii{f}_{\mathscr{L}_{\tilde{\delta}(\cdot)}(w)}.
\end{split}
\end{equation*}
\end{proof}

\section{Applications}\label{Sec5}

We will conclude this paper with some results where we essentially prove norm inequalities on the weighted variable Lebesgue spaces for some fractional type operators in the Schr\"odinger context.

First we will see how to prove that an operator $T$ is bounded in $L^{\p}(w)$ using Theorem~\ref{teo4}. These same ideas can be used to apply the other theorems.

The key point in applying Theorem 4 is to consider an appropriate $\mathcal{F}$ family. This usually requires a density argument, since we need pairs of functions $(f,g)$ such that $f$ lies both  the appropriate weighted space to apply the hypothesis and in the weighted variable Lebesgue space in which we want to obtain the thesis. 

The dense subsets of $L^{p}(w)$ are well known, for example, smooth functions and bounded functions of compact support. 
	
These sets are also dense in $L^{\p}(\R^n)$ (see for example~\cite{CUF}) and in $L^{\p}(w)$ (see~\cite{CUW}). 

More specifically, in~\cite{CUW} it is proved that if $\p\in\mathcal{P}(\R^n)$ with $p^{+}<\infty$  and $w\in L^{\p}_\loc(\R^n)$, then $L_c^{\infty}(\R^n)$, the set of bounded functions of compact support, and $C_c^\infty$, the smooth functions of compact support, are dense in $L^{\p}(w)$.

Suppose now that for all $w_0\in A_{p_0,q_0}^{\rho}$ it is verified that
\begin{equation}\label{extrapol}
\|Tf\,w_0\|_{q_0}\le C\|f\,w_0\|_{p_0}.
\end{equation}

We want to show that given a $w\in A_{\p,\q}^\rho$, $T$ maps $L^{\p}(w)$ into $L^{\q}(w)$. Since $w\in L^{\p}_\loc(\R^n)$ y $w\in L^{\q}_\loc(\R^n)$, by a standard density argument (see~\cite{CUF}, Theorem 5.39) it is sufficient to show that
\begin{equation*}
\|Tf\,w\|_{\q}\le C\|f\,w\|_{\p},
\end{equation*}
for all $f\in L_c^{\infty}$. Although intuitively, it can be thought to define $\mathcal{F}$ as
\begin{equation*}
\mathcal{F}=\{(|Tf|,|f|):\,f\in L_c^{\infty}\},
\end{equation*}
it is not known a priori that $Tf$ is in $L^{\q}(w)$. To overcome this, we again proceed by approximation and define $(Tf)_n=\min\{|Tf|,n\}\chi_{B(0,n)}$. Given that $w\in L^{\q}_{\loc}(\R^n)$, it follows that $(Tf)_n\in L^{\q}(w)$. 
On the other hand, it is clear that~\eqref{extrapol} is verified with  $|Tf|$ replaced by $(Tf)_n$. Therefore, if we define
\begin{equation*}
\mathcal{F}=\{((Tf)_n,|f|):\,f\in L_c^{\infty}, n\ge 1\},
\end{equation*}
we can apply Theorem~\ref{teo4} and Fatou's Lemma in this context (see~\cite{CUF}, Theorem 2.61) and conclude that for any $f\in L_c^{\infty}$
\begin{equation*}
\|Tf\,w\|_{\q}\le\liminf_{n\rightarrow\infty}\|(Tf)_n\,w\|_{\q}
\le C\|f\,w\|_{\p}.
\end{equation*}

\vspace{0.35cm}

\subsection{Fractional integral.} Given $0<\alpha<n$, the fractional integral of order $\alpha$ associated with $\mathcal{L}$ is defined as
\begin{equation*}\label{integralf}
\mathcal{I}_\alpha f(x) =\mathcal{L}^{-\alpha/2}f(x) 
=\int_{\R^n} K_\alpha(x,y)f(y)\,dy,
\end{equation*}
where
\begin{equation*}\label{nucleointegralf}
K_\alpha(x,y)=\int_{0}^{\infty}k_t(x,y)\,t^{\alpha/2}\,\frac{dt}{t}.
\end{equation*}
Here $k_t$ denotes the kernel of $\{e^{-t\mathcal{L}}\}_{t>0}$, the heat semigroup associated to $\mathcal{L}$.  It is known (see \cite{Kurata} and \cite{DZ-H2}) that if $V\in RH_\nu$ with $\nu>n/2$ and $\lambda_0=\min\big\{1,2-\frac{n}{\nu}\big\}$, given $N>0$ and $0<\lambda<\lambda_0$, there exists a constant $C$ such that,
\begin{equation*}\label{extradec_kt}
k_t(x,y) \leq C\,t^{-n/2}\,e^{-\frac{|x-y|^2}{C\,t}} \,\left(1+\frac{\sqrt{t}}{\rho(x)}+\frac{\sqrt{t}}{\rho(y)}\right)^{-N}\text{,}
\end{equation*}
for any $x$ and $y$ in $\R^n$. Also
\begin{equation*}\label{extradec_kt_dif}
|k_t(x,y)-k_t(x_0,y)| \leq C \left(\frac{|x-x_0|}{\sqrt{t}}\right)^{\lambda}\,
t^{-n/2}\,e^{-\frac{|x-y|^2}{C\,t}}\,\left(1+\frac{\sqrt{t}}{\rho(x)}+\frac{\sqrt{t}}{\rho(y)}\right)^{-N}\text{,}
\end{equation*}
provided that $|x-x_0|<\sqrt{t}$.

Taking into account the above estimations, the following bounding result was proved in~\cite[Proposition~8]{BCH_Lerner}.
\begin{prop}\label{prop: Lerner}
If $V\in RH_\nu$ with $\nu>n/2$, then given $N>0$ and $0<\lambda<\lambda_0$, there exists a constant $C$ such that,
\begin{equation}\label{cond_tamano_inf}
|K_\alpha(x,y)|\leq \frac{C}{|x-y|^{n-\alpha}}\bigg(1+\frac{|x-y|}{\rho(x)}\bigg)^{-N},\ \ \ x, y\in \R^n,
\end{equation}
and
\begin{equation}\label{cond_Hormander_inf}
|K_\alpha(x,y)-K_\alpha(x_0,y)|\leq C\frac{|x-x_0|^\lambda}{|x-y|^{n-\alpha+\lambda}}\left(1+\frac{|x-y|}{\rho(x_0)}\right)^{-N},
\end{equation}
whenever $|x-x_0|<\frac12|x_0-y|$.
\end{prop}

\begin{thm}\label{teo:deltaSchro}
Let $V\in RH_\nu$ with $\nu>\frac{n}{2}$, $0<\alpha<n$, $0\le\delta<\min\{\lambda_0, \alpha\}$ and $s$ such that $\frac{1}{s}=\frac{\alpha}{n}-\frac{\delta}{n}$. Then, there exists a constant $C$ such that the inequality
\begin{equation*}
\vertiii{\mathcal{I}_\alpha\,f}_{\mathscr{L}_{\delta}(w)}\le C\|f\|_{L^{s}(w^{s})},
\end{equation*}
is verified for all $w\in A^\rho_{s,\infty}$.
\end{thm}

\begin{proof}
Let $x_0\in\R^n$, and consider $\tilde{B}=B(x_0,R)$ with $R\leq\rho(x_0)$. From~\eqref{cond_tamano_inf}, the H\"older inequality and since $w^{-s'}\in A_1^{\rho}$, it follows
\begin{equation}\label{Ialpha1}
\begin{split}
\frac1{|\tilde{B}|}\int_{\tilde{B}}|\mathcal{I}_\alpha(f\chi_{2\tilde{B}})(x)|dx&\lesssim\frac1{|\tilde{B}|}\int_{\tilde{B}}\int_{2\tilde{B}}|K_\alpha(x,y)||f(y)|\,dy\,dx\\
&\lesssim\frac1{|\tilde{B}|}\int_{\tilde{B}}\int_{2\tilde{B}}\frac{|f(y)|}{|x-y|^{n-\alpha}}\bigg(1+\frac{|x-y|}{\rho(x)}\bigg)^{-N}dy\,dx\\
&\lesssim\frac{1}{|\tilde{B}|}\int_{2\tilde{B}}|f(y)|\bigg(\int_{B(y,3R)}\frac1{|x-y|^{n-\alpha}}dx\bigg)dy\\
&\lesssim |\tilde{B}|^{\alpha/n-1}\int_{2\tilde{B}}|f(y)|dy\\
&\lesssim|\tilde{B}|^{\alpha/n-1}\|fw\|_{s}\bigg(\int_{2\tilde{B}}w(y)^{-s'}dy\bigg)^{1/s'}\\
&\lesssim|\tilde{B}|^{\alpha/n-1}|\tilde{B}|^{1/s'}\frac{\|fw\|_{s}}{\|w\chi_{2\tilde{B}}\|_{\infty}}\bigg(1+\frac{2R}{\rho(x_0)}\bigg)^{\theta/s'}\\
&\leq C|\tilde{B}|^{\delta/n}\frac{\|fw\|_{s}}{\|w\chi_{\tilde{B}}\|_{\infty}}.
\end{split}
\end{equation}

We will check first condition~\eqref{promedio}. We write $f=\tilde{f}_1+\tilde{f}_2$, with $\tilde{f}_1=f\chi_{2B(x_0,\rho(x_0))}$. The observation above with $R=\rho(x_0)$ gives us the estimate we need for $\tilde{f}_1$. For $x\in B(x_0,\rho(x_0))$ and $y\in(2B(x_0,\rho(x)))^c$, we have $\rho(x)\simeq\rho(x_0)$ and $|x_0-y|\le2|x-y|$. Then, calling $B'_k=B(x_0,2^k\rho(x_0))$, by using estimate~\eqref{cond_tamano_inf} with $N>\delta+\frac{\theta}{s'}$ and H\"older inequality,  we get for $x\in B(x_0,\rho(x_0))$,
\begin{equation*}
\begin{split}
|\mathcal{I}_\alpha(\tilde{f}_2)(x)|&\lesssim\int_{(2B(x_0,\rho(x_0)))^c}\frac{|f(y)|}{|x-y|^{n-\alpha}}\bigg(1+\frac{|x-y|}{\rho(x)}\bigg)^{-N}dy\\
&\lesssim\sum_{k\geq1}\int_{B'_{k+1}\setminus B'_k}\frac{|f(y)|}{|x_0-y|^{n-\alpha}}\bigg(1+c\,\frac{|x_0-y|}{\rho(x_0)}\bigg)^{-N}dy\\
&\lesssim\frac1{\rho(x_0)^{n-\alpha}}\sum_{k\geq1}\frac1{2^{k(n-\alpha)}}\bigg(1+\frac{2^{k}\rho(x_0)}{\rho(x_0)}\bigg)^{-N}\int_{B'_{k+1}}|f(y)|\,dy\\
&\lesssim\frac{\|fw\|_{s}}{\rho(x_0)^{n-\alpha}}\sum_{k\geq1}\frac{2^{-kN}}{2^{k(n-\alpha)}}\bigg(\int_{B'_{k+1}}w(y)^{-s'}dy\bigg)^{1/s'}\\
&\lesssim\frac{\|fw\|_{s}}{\rho(x_0)^{n-\alpha}}\sum_{k\geq1}\frac{2^{-k(N-\frac{\theta}{s'})}}{2^{k(n-\alpha)}}|B'_{k+1}|^{1/s'}\left(\inf_{B'_{k+1}}w^{-s'}\right)^{1/s'}\\
&\lesssim\frac{\|fw\|_{s}}{\rho(x_0)^{n-\alpha-\frac{n}{s'}}}\sum_{k\geq1}\frac{2^{-k(N-\frac{\theta+n}{s'})}}{2^{k(n-\alpha)}}\frac{1}{\|w\chi_{B'_{k+1}}\|_{\infty}}\\
&\lesssim\frac{\|fw\|_{s}}{\|w\chi_{B(x_0,\rho(x_0))}\|_{\infty}}\rho(x_0)^{\alpha-\frac{n}{s}}\bigg(\sum_{k\geq1}2^{-k[N-(\alpha-\frac{n}{s})-\frac{\theta}{s'}]}\bigg)\\
&\lesssim\frac{\|fw\|_{s}}{\|w\chi_{B(x_0,\rho(x_0))}\|_{\infty}}|B(x_0,\rho(x_0))|^{\delta/n}.
\end{split}
\end{equation*}

Now, we must check condition~\eqref{oscilacion}. It will be enough to show that this inequality is verified for every ball $B\in\mathcal{B} _\rho$ for some constant $c_B$. Let $B=B(x_0, r)$ with $r\le\rho(x_0)$ and we write $f=f_1+f_2$, with $f_1=f\chi_{2B}$. As before, the estimate for $f_1$ follows from~\eqref{Ialpha1} for $R=r$.

Finally we will estimate $|\mathcal{I}_\alpha(f_2)(x)-c_B|$, uniformly for $x\in B$, with $c_B=\mathcal{I}_\alpha(f_2)(x_0)$. For $x\in B$ and $y\in(2B)^c$, we have $\rho(x)\simeq\rho(x_0)$ and $|x_0-y|\le2|x-y|$. Then, calling $B_k=B(x_0,2^kr)$, by using estimate~\eqref{cond_Hormander_inf} with $N=\frac{\theta}{s'}$ and $\delta<\lambda<\lambda_0$,  we have 
\begin{equation*}
\begin{split}
|\mathcal{I}_\alpha(f_2)(x)-\mathcal{I}_\alpha(f_2)(x_0)|&\lesssim\int_{(2B)^c}|f(y)|\frac{|x-x_0|^{\lambda}}{|x-y|^{n-\alpha+\lambda}}\bigg(1+\frac{|x-y|}{\rho(x)}\bigg)^{-N}dy\\
&\lesssim r^\lambda\sum_{k\geq1}\int_{B_{k+1}\setminus B_k}\frac{|f(y)|}{|x_0-y|^{n-\alpha+\lambda}}\bigg(1+c\,\frac{|x_0-y|}{\rho(x_0)}\bigg)^{-N}dy\\
&\lesssim\frac1{r^{n-\alpha}}\sum_{k\geq1}\frac1{2^{k(n-\alpha+\lambda)}}\bigg(1+\frac{2^{k}r}{\rho(x_0)}\bigg)^{-N}\int_{B_{k+1}}|f(y)|\,dy\\
&\lesssim\frac{\|fw\|_{s}}{r^{n-\alpha}}\sum_{k\geq1}\frac{1}{2^{k(n-\alpha+\lambda)}}\bigg(1+\frac{2^{k}r}{\rho(x_0)}\bigg)^{-N}\bigg(\int_{B_{k+1}}w(y)^{-s'}dy\bigg)^{1/s'}\\
&\lesssim\frac{\|fw\|_{s}}{r^{n-\alpha}}\sum_{k\geq1}\frac{|B_{k+1}|^{1/s'}}{2^{k(n-\alpha+\lambda)}}\bigg(1+\frac{2^{k}r}{\rho(x_0)}\bigg)^{-N+\frac{\theta}{s'}}\left(\inf_{B_{k+1}}w^{-s'}\right)^{1/s'}\\
&\lesssim\frac{\|fw\|_{s}}{r^{n-\alpha-\frac{n}{s'}}}\sum_{k\geq1}\frac{2^{k\frac{n}{s'}}}{2^{k(n-\alpha+\lambda)}}\frac{1}{\|w\chi_{B_{k+1}}\|_{\infty}}\\
&\lesssim\frac{\|fw\|_{s}}{\|w\chi_{B}\|_{\infty}}r^{\delta}\bigg(\sum_{k\geq1}2^{-k(\lambda-\delta)}\bigg)\\
&\lesssim\frac{\|fw\|_{s}}{\|w\chi_{B}\|_{\infty}}|B|^{\delta/n}.
\end{split}
\end{equation*}
	
\end{proof}

\begin{rem}
The above result was proved for the case $\delta=0$ in~\cite{BCH}.
\end{rem}

Thus, we can apply the Theorem~\ref{teo5} to get the following.

\begin{thm}
Let $V\in RH_\nu$ con $\nu>\frac{n}{2}$, $0<\alpha<n$ and $\p,\q\in\mathcal{P}(\R^n)$ such that
\begin{equation*}
\frac{1}{\p}-\frac{1}{\q}=\frac{\alpha}{n}.
\end{equation*}
Then, if $w\in A_{\p,\q}^{\rho}$ and  $\q\in\mathcal{P}^{\log}(\R^n)$ with $q^{-}>\frac{n}{n-\alpha}$, it follows that
\begin{equation*}
\|\mathcal{I}_\alpha\,f\|_{L^{\q}(w)}\le C\|f\|_{L^{\p}(w)}.
\end{equation*}
\end{thm}

\begin{rem}
Since $\frac{1}{\p}-\frac{1}{\q}=\frac{\alpha}{n}$ we have that  $\frac{1}{p^{+}}-\frac{1}{q^{+}}=\frac{\alpha}{n}$ from which it follows that  $p^{+}<\frac{n}{\alpha}$. The restriction $p^{+}<\frac{n}{\alpha}$ is natural for the fractional integral operator, since in the constant exponent case $\mathcal{I}_\alpha$ does not map $L^{n/\alpha}$ to $L^\infty$. 
\end{rem}

On the other hand, also as a consequence of Theorem~\ref{teo:deltaSchro}, now in combination with Theorem~\ref{teo6}, the following result follows.

\begin{thm}
Let $V\in RH_\nu$ with $\nu>\frac{n}{2}$, $0\le\delta<\min\{\alpha,\lambda_0\}$ and $s>1$ be such that $\frac1{s}=\frac{\alpha-\delta}{n}$. Let $\p\in\mathcal{P}^{\log}(\R^n)$, $p^{-}>1$, such that $\frac1{\p}=\frac{\alpha-\tilde{\delta}(\cdot)}{n}$,  $0\le\tilde{\delta}(\cdot)\textcolor{red}{<}\delta$. Then,  for every $w\in A^\rho_{\p,\infty}$ it follows that
\begin{equation*}
\vertiii{\mathcal{I}_\alpha\,f}_{\mathscr{L}_{\tilde{\delta}(\cdot)}(w)}\leq C \|f\|_{L^{\p}(w)}.
\end{equation*}
\end{thm}

\

\subsection{Commutators of Fractional Integrals} We now approach the case of the commutator of the operator $\mathcal{I}_\alpha$. For a function $b$ we will consider the commutator of  $\mathcal{I}_\alpha$ defined as,
\begin{equation*}
[\mathcal{I}_\alpha,b]f(x)=\mathcal{I}_\alpha(bf)(x)-b(x)\mathcal{I}_\alpha f (x),\ \ \ \  x\in\R^n.
\end{equation*}

It is known that for the classical case (i.e., $V\equiv0$), the commutator of the
$\mathcal{I}_\alpha$ operator is of strong type $(p,q)$ for $1<p<\frac{n}{\alpha}$ and $\frac1{p}-\frac1{q}=\frac{\alpha}{n}$ provided that $b$ belongs to $BMO$, the bounded mean oscillation space.

For the case we are dealing with here, the operator $[\mathcal{I}_\alpha,b]$, it verifies analogous bounding properties, but the class where $b$ belongs is larger than the usual $BMO$.  Let us consider the space $BMO_\rho$, defined in~\cite{BHS-Conm}, as the set of locally integrable functions $b$ such that for some $\theta>0$,
\begin{equation*}
\frac{1}{|B(x,r)|} \int_{B(x,r)} |b(y)-b_{B(x,r)}|\,dy \leq C\left(1+\frac{r}{\rho(x)}\right)^{\theta},
\end{equation*}
for all $x\in\R^n$ and $r>0$. The following result was proved in~\cite[Theorem~4.4]{Tang}. 

\begin{thm}
Let $V\in RH_\nu$ with $\nu>\frac{n}{2}$ and $b\in BMO_\rho$. Then, if $0<\alpha<n$, $1<p<\frac{n}{\alpha}$ and $\frac1{q}=\frac1{p}-\frac{\alpha}{n}$, there exists a constant $C$ such that the inequality
\begin{equation*}
\|[\mathcal{I}_\alpha,b]f\|_{L^q(w^q)}\le C\|f\|_{L^p(w^p)},
\end{equation*}
is verified for all $w\in A^\rho_{p,q}$.
\end{thm}

Hence, in view of Theorem~\ref{teo4}, we obtain the following bounding result.
\begin{thm}
Let $V\in RH_\nu$ with $\nu>\frac{n}{2}$, $b\in BMO_\rho$, $0<\alpha<n$ and  $\p,\q\in\mathcal{P}(\R^n)$ such that
\begin{equation*}
\frac{1}{\p}-\frac{1}{\q}=\frac{\alpha}{n}.
\end{equation*}
	
Then, if $w\in A_{\p,\q}^{\rho}$ and $\q\in\mathcal{P}^{\log}(\R^n)$ with $q^{-}>\frac{n}{n-\alpha}$, it follows that
\begin{equation*}
\|[\mathcal{I}_\alpha,b]f\|_{L^{\q}(w)}\le C\|f\|_{L^{\p}(w)}.
\end{equation*}
\end{thm}

\

\end{document}